
\def\LaTeX{%
  \let\Begin\begin
  \let\End\end
  \let\salta\relax
  \let\finqui\relax
  \let\futuro\relax}

\def\UK{\def\our{our}\let\sz s}
\def\USA{\def\our{or}\let\sz z}

\UK



\LaTeX

\UK


\salta

\documentclass[twoside,12pt]{article}
\setlength{\textheight}{24cm}
\setlength{\textwidth}{16cm}
\setlength{\oddsidemargin}{2mm}
\setlength{\evensidemargin}{2mm}
\setlength{\topmargin}{-15mm}
\parskip2mm


\usepackage[usenames,dvipsnames]{color}
\usepackage{amsmath}
\usepackage{amsthm}
\usepackage{amssymb, bbm}
\usepackage[mathcal]{euscript}

\usepackage{hyperref}
\usepackage{enumitem}
\usepackage{graphicx}

%
%


\definecolor{viola}{rgb}{0.3,0,0.7}
\definecolor{ciclamino}{rgb}{0.5,0,0.5}
\definecolor{rosso}{rgb}{0.85,0,0}

\def\anold #1{{\color{rosso}#1}}

\def\anold #1{#1}




\bibliographystyle{plain}


%
\newtheorem{theorem}{Theorem}
\newtheorem{remark}{Remark}

\finqui

\def\Bcenter{\Begin{center}}
\def\Ecenter{\End{center}}
\let\non\nonumber




\def\step #1 \par{\medskip\noindent{\bf #1.}\quad}


\def\Holder{H\"older}

\def\aand{\quad\hbox{and}\quad}

\def\rhs{right-hand side}



\def\multibold #1{\def\arg{#1}%
  \ifx\arg\pto \let\next\relax
  \else
  \def\next{\expandafter
    \def\csname #1#1#1\endcsname{{\boldsymbol #1}}%
    \multibold}%
  \fi \next}

\def\pto{.}

\def\multical #1{\def\arg{#1}%
  \ifx\arg\pto \let\next\relax
  \else
  \def\next{\expandafter
    \def\csname #1#1\endcsname{{\cal #1}}%
    \multical}%
  \fi \next}


\def\multimathop #1 {\def\arg{#1}%
  \ifx\arg\pto \let\next\relax
  \else
  \def\next{\expandafter
    \def\csname #1\endcsname{\mathop{\rm #1}\nolimits}%
    \multimathop}%
  \fi \next}

\multibold
qwertyuiopasdfghjklzxcvbnmQWERTYUIOPASDFGHJKLZXCVBNM.

\multical
QWERTYUIOPASDFGHJKLZXCVBNM.


\multimathop
diag dist div dom mean meas sign supp .


\def\Accorpa #1#2 #3 {\gdef #1{\eqref{#2}-\eqref{#3}}%
  \wlog{}\wlog{\string #1 -> #2 - #3}\wlog{}}
\def\Accorparef #1#2 #3 {\gdef #1{\ref{#2}--\ref{#3}}%
  \wlog{}\wlog{\string #1 -> #2 - #3}\wlog{}}


\def\<#1>{\mathopen\langle #1\mathclose\rangle}
\def\norma #1{\mathopen \| #1\mathclose \|}

\def\[#1]{\mathopen\langle\!\langle #1\mathclose\rangle\!\rangle}

\def\iot {\int_0^t}

\def\iO{\int_\Omega}

\def\dt{\partial_t}
\def\dn{\partial_\nnn}

\def\checkmmode #1{\relax\ifmmode\hbox{#1}\else{#1}\fi}
\def\aeO{\checkmmode{a.e.\ in~$\Omega$}}

\def\aet{\checkmmode{a.e.\ in~$(0,T)$}}


\def\erre{{\mathbb{R}}}
\def\enne{{\mathbb{N}}}




\def\genspazio #1#2#3#4#5{#1^{#2}(#5,#4;#3)}
\def\spazio #1#2#3{\genspazio {#1}{#2}{#3}T0}

\def\L {\spazio L}
\def\H {\spazio H}

\def\C #1#2{C^{#1}([0,T];#2)}


\def\Lx #1{L^{#1}(\Omega)}
\def\Hx #1{H^{#1}(\Omega)}

\def\Hdn{{H^2_\nnn(\Omega)}}



\let\theta\vartheta
\let\badeps\epsilon
\let\eps\varepsilon
\let\ph\varphi

\let\TeXchi\chi                         
\newbox\chibox
\setbox0 \hbox{\mathsurround0pt $\TeXchi$}
\setbox\chibox \hbox{\raise\dp0 \box 0 }
\def\chi{\copy\chibox}



\def\phz{\ph_0}
\def\sigmaz{\sigma_0}

\def\soluz{(f, \ph,\mu,\sigma)}

\def\Vp{(\Hx1)^*}

\def\MM{{\cal M}}

\def\mobm{\anold{{\mathbbm{m}}}}

\def\hh{{\mathbbm{h}}}
\def\pp{\anold{{\mathbbm{p}}}}
\def\qq{\anold{{\mathbbm{q}}}}
\def\kk{\anold{{\mathbbm{k}}}}
\def\ww{\anold{{\mathbbm{w}}}}
\def\ov #1{{\overline{#1}}}

\def\bv{{\boldsymbol{v}}}

\def\Rcel{\RR}
\def\Rold{{{H}}}
\Begin{document}


%
\title{On a phenotype-structured phase-field model of nutrient-limited tumour growth}
\date{}
\author{}
\maketitle
\Bcenter
\vskip-1cm
{\large\sc Tommaso Lorenzi$^{(1)}$}\\
{\normalsize e-mail: {\tt tommaso.lorenzi@polito.it}}\\[.25cm]
{\large\sc	Giulia Pozzi$^{(1)}$ }\\
{\normalsize e-mail: {\tt giulia.pozzi@polito.it}}\\[.25cm]
{\large\sc Andrea Signori$^{(2)}$}\\
{\normalsize e-mail: {\tt andrea.signori@polimi.it}}\\[.5cm]
$^{(1)}$
{\small Department of Mathematical Sciences \lq\lq G. L. Lagrange"}\\
{\small Politecnico di Torino, 10129 Torino, Italy}\\[.3cm] 
$^{({2})}$
{\small Dipartimento di Matematica, Politecnico di Milano}\\
{\small via Bonardi 9, 20133 Milano, Italy}
\\
{Alexander von Humboldt Research Fellow}
%
%
%
\date{}

\Ecenter

\Begin{abstract}
\noindent
Phase-field models of tumour growth have proved useful as theoretical tools to investigate cancer invasion. A key implicit assumption underlying mathematical models of this type which have so far been proposed, though, is that cells in the tumour are identical. This assumption ignores both the fact that cells in the same tumour may express different characteristics to different extents, exhibiting heterogeneous phenotypes, and the fact that cells may undergo phenotypic changes, with their characteristics evolving over time. To address such a limitation, in this paper we incorporate inter-cellular phenotypic heterogeneity and the evolution of cell phenotypes into the phase-field modelling framework. This is achieved by formulating a phenotype-structured phase-field model of nutrient-limited tumour growth. For this model, we first establish a well-posedness result under general assumptions on the model functions, which encompass a wide spectrum of biologically relevant scenarios. We then present a sample of numerical solutions to showcase key features of spatiotemporal and evolutionary cell dynamics predicted by the model. We conclude with a brief overview of modelling and analytical research perspectives.

\End{abstract}
\salta
\pagestyle{myheadings}
\newcommand\testopari{\sc Lorenzi \ -- \ Pozzi \ -- \ Signori }
\newcommand\testodispari{{\sc On a phenotype-structured phase-field model of tumour growth}}
\markboth{\testodispari}{\testopari}
\finqui
%

\section{Introduction}
\label{SEC:INTRO}
\setcounter{equation}{0}

\subsection{Background}
Phase-field models of tumour growth have been increasingly used to investigate different aspects of cancer dynamics. Originally proposed for studying inert matter processes in the theoretical framework of out-of-equilibrium thermodynamics, phase-field models have then proved useful also in studying morphogenesis and growth processes in living matter, as a result of their ability to describe nonlinear spatiotemporal phenomena in the presence of moving boundaries~\cite{travasso2011phase}. A variety of mathematical models of this type, based on different extensions of the Cahn--Hilliard equation~\cite{CH1,CH2}, have been proposed to capture different aspects of tumour growth and therapy. Amongst these, we mention models incorporating chemotaxis and active transport \cite{GLSS,RSchS}, tumour-microenvironment interactions~\cite{chatelain2011morphological,chen2014tumor, cristini2009nonlinear, wise2008three}, chemotherapy and radiotherapy \cite{agosti2018personalized}, and immunotherapy \cite{pozzi2022t}.

As a prototypical example of phase-field models of tumour growth, we consider the following system modelling nutrient-limited growth~\cite{GLSS}:
\begin{alignat}{2}
	\label{eq:01}
	& \dt \ph - \div \left(\mobm(\ph) \nabla \mu\right)
	=
	\hh(\ph) \left( \sigma \, p - q \right)
	\qquad&& \text{in $Q$,} \phantom{\int_\YY}
	\\
	\label{eq:02}
	& \mu = - \eps \Delta \ph + \dfrac 1\eps F' (\ph)
	\qquad&& \text{in $Q$,} 
	\\
	\label{eq:04}
	& \dt \sigma - D_{\sigma} \, \Delta \sigma
	+ \hh(\ph) \, k \, \sigma
	=
	b\left(\sigma_B - \sigma\right)
	\qquad&& \text{in $Q$,} \phantom{\int_\YY}
\end{alignat}
subject to the zero-flux boundary conditions 
\begin{alignat}{2}	
	\label{eq:05}
	& \dn \ph=  (\mobm(\ph) \nabla \mu) \cdot \nnn = \dn \sigma
	=
	0
	\qquad && \text{on $\Sigma$},
\end{alignat}
and to the initial conditions
\begin{alignat}{2}
	\label{eq:06}
	& \ph(0)=\ph_0, \quad \sigma(0)=\sigma_0
	\qquad && \text{in $\Omega$.}
\end{alignat}
Here $\Omega\subset\erre^d$, with $d \in \{2,3\}$, is a bounded spatial domain with smooth boundary $\Gamma:=\partial\Omega$, $\nnn$ indicates the outward unit normal on $\Gamma$, with corresponding directional derivative $\dn$, and, for a fixed final time $T \in \mathbb{R}^+$, $Q$ and $\Sigma$ denote, respectively, the parabolic cylinder and its boundary, that is, 
\begin{align*}
  Q := \Omega\times(0,T),
  \qquad
  \Sigma := \Gamma\times(0,T).
\end{align*}

In the phase-field model~\eqref{eq:01}-\eqref{eq:04}, the function $\ph =\ph(x,t)$ acts as an order parameter, indicating the difference in volume fractions between tumour cells and healthy cells at position $x \in \Omega$ and time $t \in (0,T)$. In particular, it is common to consider $\ph: Q \to [-1,1]$, with: $\ph(x,t)=-1$ if only healthy cells are present at $(x,t)$; $\ph(x,t)=1$ if there are only tumour cells at $(x,t)$; and $-1<\ph(x,t)<1$ if there is a mixture of healthy and tumour cells at $(x,t)$, with values of $\ph(x,t)$ closer to $1$ corresponding to a higher fraction of tumour cells. Moreover, as it is standard in the framework of the Cahn--Hilliard equation, the function $\mu =\mu(x,t)$ represents the chemical potential associated with the phase variable $\ph$. Finally, the non-negative function $\sigma =\sigma(x,t)$ models the local concentration of nutrients (e.g. oxygen and glucose), which are consumed by the cells to produce the energy required by cell division.

In~\eqref{eq:01}, the parameter $p \in \mathbb{R}^+$  models the rate of cell proliferation, which is multiplied by $\sigma$ to take into account the fact that cell proliferation is dependent
on availability of nutrients which fuel cell division (i.e. arrest of tumour growth occurs due to nutrient insufficiency), and the parameter $q \in \mathbb{R}^+$ is the cell death rate. Moreover, the second term on the left-hand side of~\eqref{eq:01}, where $\mobm(\ph)$ is a non-negative function linked to cell mobility, represents the rate of change of $\ph$ due to cell movement. 

In~\eqref{eq:04}, the parameter $D_{\sigma} \in \mathbb{R}^+$  is the diffusion coefficient of the nutrients, while the parameter $k \in \mathbb{R}^+$ is the rate of consumption of the nutrients by the cells. Furthermore, the term on the right-hand side of~\eqref{eq:04}, where $b \in \mathbb{R}^+$ and $\sigma_B$ is a non-negative function, represents the rate of change of $\sigma$ due to the supply and decay of the nutrients. 

In~\eqref{eq:01} and~\eqref{eq:04}, the function $\hh(\ph)$ is a continuous and bounded truncation function, on which the following assumptions are commonly made
\begin{align}
	\label{ass:h}
\hh(-1) = 0, \quad \hh(1) = 1, \quad 0 < \hh(\ph) \leq 1 \; \text{ for } \; \ph \in (-1,1),
\end{align}
with a prototypical example being $\hh(\ph) := \frac 12 (\ph+1).$
In view of assumptions~\eqref{ass:h}, this function is used in~\eqref{eq:01} and~\eqref{eq:04} to translate into mathematical terms the ideas that the effects of proliferation and death of healthy cells on the dynamic of $\ph$ and the effect of nutrient consumption by healthy cells on the dynamic of $\sigma$ can be, as a first approximation, neglected. This modelling choice is rooted in the idea that proliferation and death of healthy cells can be considered as occurring on a slower time scale compared to proliferation and death of tumour cells, and the rate of nutrient consumption of tumour cells can be assumed to be much larger than that of healthy cells.

Finally, the parameter $\eps \in (0,1)$ in~\eqref{eq:02} is linked to the thickness of the transition layer between the two pure phases (i.e. $\ph=\pm 1$) and the term $F'$ denotes the derivative of a potential $F$, which exhibits a double-well shape. Examples of such potentials include regular potentials, such as the standard quartic regular potential, $F_{\rm reg}$, and the logarithmic potential, $F_{\rm log}$, which are defined, respectively, as
\begin{align}
	\label{regpot}
	& F_{\rm reg}(r) := \frac{1}{4}(r^2-1)^2, 
	\quad r \in \mathbb{R},
	\\ \label{Flog}
	& F_{\rm log}(r)
	:= (1+r)\ln(1+r)+(1-r)\ln (1-r) - c_1 r^2, \quad  r \in (-1,1),
\end{align}
where $c_1 \in \mathbb{R}^+$ with $c_1>1$. Singular potentials are often approximated by regular ones, inducing slight deviations of the values of $\ph$ from the physical range $[-1,1]$. A classical example is the regular quartic potential in~\eqref{regpot}. When regular potentials are considered, the physical property $\ph \in [-1,1]$ does not hold rigorously, but it is satisfied only approximately. From a numerical standpoint, however, such a property can be maintained with high accuracy by choosing a sufficiently small interface thickness, which is directly linked to the choice of the parameter $\varepsilon$.

\subsection{Object of study}
A key implicit assumption underlying phase-field models of tumour growth which have so far been proposed, such as~\eqref{eq:01}-\eqref{eq:04}, is that cells in the tumour are identical -- namely, they exhibit the same proliferation and death rates and identical nutrient consumption. This assumption ignores both the fact that cells in the same tumour may express different characteristics to different extents, meaning they have different phenotypes, and the fact that cells may undergo phenotypic changes, and thus their characteristics evolve over time. Since these biological facts are known to be crucial in cancer development and progression~\cite{Marusyketal2012}, neglecting them limits the applicability of mathematical models that rely on this simplifying assumption.

A possible approach to incorporate inter-cellular phenotypic heterogeneity and the evolution of cell phenotypes into continuum models of tumour development consists in: representing the cell phenotypic state by means of a continuous structuring variable, which captures variability in certain phenotypic characteristics amongst the cells; and conceptualising phenotypic changes as transitions between phenotypic states, which are modelled through an integral or a differential operator -- see the review~\cite{Chisholmetal2016} and references therein. In this modelling framework, the evolution of the cell distribution in the space of phenotypic states (i.e. the cell phenotype distribution) is governed by an integro-differential equation or a partial differential equation (see the review~\cite{Lorenzietal2025} and references therein), which provides a mathematical depiction of the cell evolutionary dynamics by describing how the cell phenotype distribution is dynamically shaped by the interplay between environmental selection, acting upon the phenotypic expression of the cells, and by phenotypic changes undergone by the cells.

Building on this approach, in order to address the aforementioned limitation of existing phase-field models of tumour growth, we propose the following phenotype-structured version of the phase-field model~\eqref{eq:01}-\eqref{eq:04}:
\begin{alignat}{2}
	\label{eq:1}
	& \dt \ph - \div \left(\mobm(\ph) \nabla \mu\right)
=
	\hh(\ph) \left( \sigma \int_\YY \pp f \,{{\rm d}y} 
	-\int_\YY \qq f \,{{\rm d}y} \right)
	\qquad && \text{in $Q$,} \phantom{\int_\YY}
	\\
	\label{eq:2}
	& \mu = - \eps \Delta \ph + \dfrac 1\eps F' (\ph)
	\qquad && \text{in $Q$,} \phantom{\int_\YY}
	\\
	\label{eq:3}
	& \dt f = 
 \alpha \, \hh (\ph) \, \left[\theta \, \left(\int_{\cal Y} {\MM}f \,{{\rm d}y} - f\right) + \left(\Rcel - \int_\YY \Rcel f \,{\rm d}y\right) f\right]	\qquad && \text{in $Q \times \YY$,} \phantom{\int_\YY}
	\\
	\label{eq:4}
	& \dt \sigma - D_{\sigma} \, \Delta \sigma
+ \hh(\ph) \, \sigma \, \left(\int_\YY \kk f \,{{\rm d}y}\right)
	=
	b\left(\sigma_B - \sigma\right)
	\qquad && \text{in $Q$,}
\end{alignat}	
which we complement with the following zero-flux boundary conditions 
\begin{alignat}{2}	
	\label{eq:5}
	& \dn \ph=  (\mobm(\ph) \nabla \mu) \cdot \nnn = \dn \sigma
	=
	0
	\qquad && \text{on $\Sigma$}{,}
	\end{alignat}
and the initial conditions	
\begin{alignat}{2}	
	\label{eq:6}
	& \ph(0)=\ph_0, \quad \sigma(0)=\sigma_0
	\qquad && \text{in $\Omega$,} \phantom{\int_\YY}
	\\
	\label{eq:7}
	& f(0)=f_0
	\qquad && \text{in $\Omega \times \YY$.} 	
\end{alignat}
\Accorpa\Sys {eq:1} {eq:7}
Compared to the phase-field model~\eqref{eq:01}-\eqref{eq:04}, in the phenotype-structured phase-field model~\eqref{eq:1}-\eqref{eq:4} there are: 
\begin{itemize}
\item[(i)] An additional independent variable $y \in \YY \subseteq \mathbb{R}$, which is a continuous structuring variable representing the cell phenotype and capturing inter-cellular variability in proliferation and death rates and in the rate of consumption of nutrients. Thus, $\YY$ denotes the phenotype domain.
\item[(ii)] An additional dependent variable $f =f(x,t,y)$, which is a probability density function, i.e. 
\begin{align}
	\label{proprel:f}
f \geq 0, \qquad \int_{\cal Y} f(\cdot,\cdot,y) \,{{\rm d}y} = 1,
\end{align}
that represents the distribution of cells over $\YY$ at position $x \in \Omega$ and time $t \in (0,T)$ (i.e. the local cell phenotype distribution). 
\end{itemize}
Moreover, the functions $\pp=\pp(y)$ and $\qq=\qq(y)$ in~\eqref{eq:1} and the function $\kk=\kk(y)$ in~\eqref{eq:3}, which replace, respectively, the parameters $p$ and $q$ in~\eqref{eq:01} and the parameter $k$ in~\eqref{eq:04}, are non-negative phenotype-dependent functions that model, in analogy with the corresponding parameters, the rate of proliferation, the death rate, and the rate of consumption of nutrients of cells in the phenotypic state $y$. Hence, the integral terms  
$$
\int_\YY \pp(y) f(x,t,y) \,{{\rm d}y}, \qquad \int_\YY \qq(y) f(x,t,y) \,{{\rm d}y}, \qquad \int_\YY \kk(y) f(x,t,y) \,{{\rm d}y},
$$
in~\eqref{eq:1} and~\eqref{eq:3} represent the mean proliferation rate, the mean death rate, and the mean rate of consumption of nutrients of cells at position $x$ and time $t$, respectively. 

In the first term on the right-hand side of~\eqref{eq:3}, the parameter $\theta \in \mathbb{R}^+$ models the rate at which cells undergo phenotypic changes and $\MM =\MM(x, y', y)$ is a probability kernel, i.e. 
$$
\MM \geq 0, \qquad \int_{\cal Y} \MM(\cdot,\cdot,y) \,{{\rm d}y} = 1,
$$  
that models the probability with which phenotypic changes lead cells at position $x$ to transition from the phenotypic state $y'$ to the phenotypic state $y$. 
The dependence on the variable $x$ of the function $\MM$ takes into account the fact that the way in which phenotypic changes take place may vary across the spatial domain $\Omega$ due to spatial variability in the local microenvironment to which the cells are exposed~\cite{Huang2013}. Moreover, in the vein of~\cite{Hameletal2020,TsimringLevine1996}, in the second term on the right-hand side of~\eqref{eq:3}, the function $\RR =\RR(x,t,y)$ models the fitness of cells in the phenotypic state $y$ under the environmental conditions of position $x$ at time $t$, that is, the cell fitness landscape~\cite{Burger2000,lorenzipouchol2020,villaetal2021}, while the integral term 
$$
\int_\YY \Rcel(x,t,z) f(x,t,z) \,{{\rm d}z}
$$
represents the mean fitness of cells at position $x$ and time $t$. The second term on the right-hand side of~\eqref{eq:3} relies on the observation that if the fitness of a phenotype $y$ is smaller/larger than the average fitness at $(x,t)$, then the frequency of the phenotype $y$ decreases/increases at $(x,t)$. As similarly done in~\eqref{eq:1} and~\eqref{eq:3}, in the light of assumptions~\eqref{ass:h}, the function $\hh(\ph)$ is introduced in the right-hand side of~\eqref{eq:3} to translate into mathematical terms the idea that, since the evolutionary rate of healthy cells can be assumed to be negligible compared to that of tumour cells, the effect of the evolutionary dynamics of healthy cells on the dynamic of $f$ can be, as a first approximation, neglected. Finally, the parameter $\alpha \in \mathbb{R}^+$ measures the timescale over which the phenotypic evolution of tumour cells occur.

We notice that, in the framework of the model~\eqref{eq:1}-\eqref{eq:4}, the local mean phenotype and the local phenotype variance (i.e. the mean phenotype of cells at position $x$ and the corresponding variance) at time $t$ are defined, respectively, as
\begin{equation}
\label{def:meanvarphen}
\int_\YY y \, f(x,t,y) \,{\rm d}y \quad \text{and} \quad \int_\YY y^2 \, f(x,t,y) \,{\rm d}y - \left(\int_\YY y \, f(x,t,y) \,{\rm d}y \right)^2.
\end{equation}

\begin{remark}
Note that, as it is shown in the paragraph {\bf Energy estimate} in Section~\ref{SEC:EX}, a dissipation inequality analogous to the one derived in~\cite{GLSS} for the phase-field model~\eqref{eq:01}-\eqref{eq:04} can also be derived for the phenotype-structured phase-field model~\eqref{eq:1}-\eqref{eq:4}, thus indicating that this model fulfils the second law of thermodynamics.
\end{remark}

\subsection{Content of the paper}
In this paper, we study the well-posedness of the initial-boundary value problem defined by complementing the phenotype-structured phase-field model of nutrient-limited tumour growth~\eqref{eq:1}-\eqref{eq:4} with the boundary conditions~\eqref{eq:5} and the initial conditions~\eqref{eq:6}-\eqref{eq:7}. Moreover, as a proof of concept for the ideas underlying this modelling framework, we present a sample of numerical solutions, which showcase key features of spatiotemporal and evolutionary cell dynamics predicted by the model.

The remainder of the paper is organised as follows. In Section~\ref{SEC:NOT}, we introduce a few technical preliminaries and the main assumptions on the model functions. In Section~\ref{SEC:MARE}, Theorem~\ref{THM:WP} establishes the existence of weak solutions, with the proof given in Subsection~\ref{SEC:EX}, while Theorem~\ref{THM:UQ} establishes uniqueness, with its proof provided in Subsection~\ref{SEC:UQ}.
In Section~\ref{SEC:NUM}, we first describe the scheme, which relies on the finite element method, employed to solve numerically the initial-boundary value problem of the model (see Subsection~\ref{SUBSEC:methods}) along with the simulation set-up (see Subsection~\ref{SUBSEC:setup}), before presenting the main results of numerical simulations (see Subsection \ref{SUBSEC:results}). Section~\ref{SEC:CONC} concludes the paper and provides a brief overview of modelling and analytical research perspectives.

\vskip3mm

\section{Notation, preliminaries, and main assumptions}
\label{SEC:NOT}

\subsection{Notation and preliminaries}
{We start by recalling that, as mentioned earlier, here $\Omega \subset \mathbb{R}^d$, with $d \in \{2,3\}$, is a bounded spatial domain with a smooth boundary $\Gamma := \partial\Omega$, and we use the notation $|\Omega|$ to denote the Lebesgue measure of $\Omega$.

For any Banach space $X$, we denote by $\norma{\cdot}_X$ the norm of $X$, by $X^*$ the dual space of $X$, and by $\<\cdot,\cdot>_X$ the duality pairing between $X^*$ and $X$. Moreover, when $X$ is a Hilbert space, we write the inner product as $(\cdot,\cdot)_X$.
The standard Lebesgue and Sobolev spaces defined on $\Omega$ for each $1 \leq p \leq \infty$ and $k \geq 0$ are written as $L^p(\Omega)$ and $W^{k,p}(\Omega)$, with their respective norms $\norma{\cdot}_{L^p(\Omega)}=:\norma{\cdot}_{p}$ and $\norma{\cdot}_{W^{k,p}(\Omega)}$. For the sake of brevity, the norm of $\Lx2$ is simply denoted as $\norma{\cdot}$ and the duality pairing of $\Hx1$ is denoted by $\<\cdot,\cdot>$.
Similar notations are employed for Lebesgue and Sobolev spaces defined on $\Gamma$, $Q$, $\Sigma$ or $\YY$.

Moreover, we adopt the standard convention $H^k(\Omega):= W^{k,2}(\Omega)$ for all $k\in\mathbb{N}$, and we set
\begin{align*}
  \Hdn:=\{ v\in H^2(\Omega): \dn v=0\;\,\text{a.e.~on } \Gamma \}.
\end{align*}
Furthermore, for a generic element $v \in \Vp$, we define its generalised mean value as $v_\Omega~:=~|\Omega|^{-1} \, \<v,1>$.

Finally, we introduce a mathematical tool that is commonly employed in the study of problems related to Cahn--Hilliard-like equations. Consider the following problem: for a given $\psi \in \Vp$, seek $u \in \Hx1$ such that
\begin{align}
	\label{weak:neu}
	\iO \nabla u \cdot \nabla v = \< \psi , v > 
	\quad \forall v \in \Hx1.
\end{align}
This corresponds to the usual weak formulation of the homogeneous Neumann problem for the Poisson equation $-\Delta u = \psi$ when $\psi \in \Lx2$.
The solvability of \eqref{weak:neu} for $\psi~\in~\Vp$ is contingent on the condition 
that $\psi$ possesses a zero mean value. Furthermore, if this condition is met, a unique solution with zero mean value exists. Consequently, the operator 
${\cal N}: \{\psi \in \Vp : \psi_\Omega = 0\} \to \{u \in \Hx1 : u_\Omega = 0\}$, defined by mapping $\psi$ onto the unique solution $u$ to \eqref{weak:neu} with $u_\Omega = 0$, is well-defined. This operator establishes an isomorphism between the aforementioned function spaces. In addition, the norm $\norma{\psi}_{\Vp}^2 := \norma{\nabla {\cal N} (\psi - {\psi_\Omega})}^2 + |{\psi_\Omega}|^2$ for $\psi \in \Vp$ can be shown to define a Hilbert norm in $\Vp$ equivalent to the standard dual norm. From these definitions, it follows that
\begin{align*}
  &\iO \nabla {\cal N}\psi \cdot \nabla v = \< \psi , v > \quad \text{for every $\psi \in \dom({\cal N})$ and $v \in \Hx1$},\\
  &\< \psi , {\cal N}\zeta > = \< \zeta , {\cal N}\psi > \quad \text{for every $\psi, \zeta \in \dom({\cal N})$}, \\
  &\< \psi , {\cal N}\psi > = \iO |\nabla{\cal N}\psi|^2 = \norma{\psi}_{\Vp}^2 \quad \text{for every $\psi \in \dom({\cal N})$}.
\end{align*}
Moreover, it also holds that
\begin{align*}
	\iot \< \dt v(s) , {\cal N} v(s) > \, ds 
	= \iot \< v(s) , {\cal N}(\dt v(s)) > \, ds 
	= \frac 12 \, \norma{v(t)}_{\Vp}^2 - \frac 12 \, \norma{v(0)}_{\Vp}^2 
\end{align*}
for every $t \in [0, T]$ and every $v \in \H1 {\Vp}$ satisfying $v_\Omega = 0$ almost everywhere.

\subsection{Main assumptions}
We make the following assumptions {\rm \bf A1}-{\rm \bf A7} on the model functions, which are general enough to encompass a wide range of biological scenarios.

\begin{enumerate}[label={\bf A\arabic{*}}, ref={\bf A\arabic{*}}]
\item  \label{ass:0:mob}
The function $\mobm \in C^0(\erre)$ is globally Lipschitz continuous and there exist some positive real constants $m^*$ and $M^*$ such that 
\begin{align}
	\label{mob:pos}
	0 <m^*\leq \mobm(r) \leq M^* <\infty
		\quad \text{for all $r\in\erre$.} 
\end{align}

\item  \label{ass:1:potenziale}
The potential $F:\erre\to[0,+\infty)$ is twice differentiable and can be decomposed as $F= F_1+F_2$, with $F_1$ nonnegative, convex and $F_2 \in C^2(\erre)$ a quadratic perturbation.
Namely, it holds that, for some $c_0\in \erre^+$,
\begin{align*}
	|F_2(r)| \leq c_0(r^2 +1 ),
	 \quad r \in \erre,	
	\qquad 
	|F''_2(r)| \leq c, 
	 \quad r \in \erre,	
\end{align*}
and the following conditions hold for all $r\in\erre$:
    \begin{alignat}{2}
    \label{growth:F:1}
       &\exists \, c_1 \in \mathbb{R}^+ : \,\,  {F(r)} \geq 
        c_1(1+r^2),
         \\
        \label{growth:F:2}
     &\forall \badeps \in \mathbb{R}^+ \,\, \exists \, c_\badeps \in \mathbb{R} \, :  \,\,|{F'(r)}| \le
        c_{\badeps} + \badeps F(r).
    \end{alignat}
    Without loss of generality, we scale $F_1$ by assuming $F_1(0)=F_1'(0)=0.$

\item  \label{ass:2:nuclei}
The probability kernel $\MM$, with density also denoted by $\MM$, satisfies the following assumptions
\begin{align}
	& \label{prop:M:1}
	\mathcal{M} : \Omega \times \YY \to {\rm P}(\YY), 
	\quad 
	\text{with } \mathcal{M}(x, y)(dy') = \mathcal{M}( x, y,y')\,dy', \\
	& \label{prop:M:2}
	\text{ where } \mathcal{M} : \Omega \times \YY \times \YY \to [0,\infty) \text{ is measurable,} \\
	& \label{prop:M:3}
	\text{and } \int_\YY \mathcal{M}(x, y,y')\,dy'= 1 \quad \text{for a.e. } (x,y) \in \Omega \times \YY.
\end{align}
Here, $ {\rm P}(\YY) $ denotes the space of Borel probability measures on $ \YY $. Moreover, we assume $\theta\in \erre^+$.

\item  \label{ass:3:pheno}
The functions $\pp$, $\qq$, $\kk$, and $\ww$ satisfy the following assumptions 
\begin{align}
	& \label{ass:functions}
	\pp,\qq,\kk,\ww: \YY \to \erre, \quad
	\pp,\qq,\kk,\ww \in L^\infty (\YY) \cap L^1(\YY).
\end{align}

\item  \label{ass:new:erre}
The function $\Rcel$ satisfies the following assumptions 
\begin{align}
	&  \label{prop:R:1}
	{\Rcel}: Q \times \YY \to \erre,
	\quad 
	{\Rcel} \in L^\infty(Q) \times (L^\infty(\YY) \cap L^1(\YY) ).
\end{align}

\item  \label{ass:4:troncante}
The function $\hh\anold{:\erre \to \erre}$ is nonnegative, bounded, and continuous. 

\item \label{ass:5:sB}
	The prescribed function $\sigma_B$ is such that $\sigma_B\in \L\infty{\Lx2}$.

\end{enumerate}
\Accorparef\tutte {ass:0:mob} {ass:5:sB}

\begin{remark}
Note that, while the quartic potential defined via~\eqref{regpot} does fulfil assumptions~\ref{ass:1:potenziale}, the same assumptions are also satisfied by higher-order polynomial-type potentials and exponential potentials, but are not satisfied by singular  potentials such as~\eqref{Flog}. Note also that, for the existence of weak solutions, the explicit structure of $\hh$ as a truncating function is not necessary, and more general forms can be considered.
\end{remark}%

\section{A well-posedness result}
\label{SEC:MARE}
Without loss of generality, throughout this section we set the model parameters $\eps$, $D_{\sigma}$, and $\alpha$ to unity for convenience.

Within the framework of assumptions \tutte, the existence of weak solutions to the initial-boundary value problem~\eqref{eq:1}-\eqref{eq:7} is established by the following theorem.

\begin{theorem}
\label{THM:WP}
Let assumptions \tutte\ hold and assume that the initial data satisfy the following assumptions
\begin{align}
	\label{ass:data}
	& \ph_0 \in {\Hx1} ,
	\quad 
	F(\ph_0) \in \Lx1,	
	\quad 
	\sigma_0 \in \Hx1 \cap \Lx\infty,
	\\
	\label{f0:density:1}
	& 
	\quad 
	{f_0(\cdot,y)  \geq 0 \,\, \text{for a.e. $y \in \YY$},}
	\quad f_0 \in L^\infty(\Omega)\times {\big(}L^\infty(\YY) \cap L^1(\YY){\big)}, 
	\\ &  
	 \int_\YY f_0(x,y)\,{{\rm d}y} =1 \quad \text{for a.e. $x \in \Omega.$}
	\label{f0:density:2}
\end{align}
Then there exists a weak solution to the initial-boundary value problem \Sys, in the following sense. There exists a quadruplet $\soluz$ which enjoys the regularity properties
\begin{align}
  \label{reg:0:f}
  & f \in \C0 {L^\infty(\Omega)\times L^1(\YY) },
  \quad 
   {f(\cdot,\cdot,y)\geq 0\,\,  \text{for a.e. $y \in \YY$,}  }
  \\ 
  & 
  \int_\YY f(x,t,y){\, {{\rm d}y}} =1 \,\, \text{for all $t \in [0,T]$ and a.e. $x \in \Omega$},
  \label{reg:0:f:nonneg}
  \\
  & \ph \in \H1 \Vp \cap \L\infty {\Hx1} \cap \L2 {\Hdn}, 
  \label{reg:1:phi}
  \\
  & \mu \in \L2 {\Hx1}, 
  \label{reg:2:mu}
  \\
  & \sigma \in \H1 {\Lx2} \cap \L\infty {\Hx1} \cap \L2 {\Hdn} \cap L^\infty(Q),
  \label{reg:3:sigma}
\end{align}
and solves the variational identities
\begin{align}
  & \< \dt \ph , v>
  + \iO \mobm(\ph) \nabla  \mu \cdot \nabla v
  = \iO \hh(\ph) \left( \sigma \int_\YY \pp f \,{{\rm d}y} 
	-\int_\YY \qq f \,{{\rm d}y} \right) v 
  \non
  \\
  & \quad \hbox{for every $v\in {\Hx1}$ and \aet},
  \label{wf:1}
  \\
  & 
  \iO \mu v
  = 
  {\iO} \nabla\ph \cdot \nabla v
  + \iO F'(\ph) v
  \non
  \\
  & 
  \quad \hbox{for every $v\in {\Hx1}$ and \aet},
  \label{wf:2}
  \\
  & \<\dt\sigma , v>
  + \iO \nabla\sigma \cdot \nabla v
  + \iO \hh(\ph)\sigma \left(\int_\YY \kk f \,{{\rm d}y} \right)v
  = b \iO \bigl(\sigma_B - \sigma\bigr) v
  \non
  \\
  & 
  \quad \hbox{for every $v\in {\Hx1}$ and \aet},
  \label{wf:3}
  \\
  & \ph(0) = \phz
  \aand
  \sigma(0) = \sigmaz
  \label{wf:4}
  \quad \aeO,
\end{align}
with 
{
\begin{align}
  \non 
  f(t) & = f_0 + \theta\iot \hh (\ph(s))  \left( \int_{\cal Y} {\MM}(s,y')f(s,y') \,{{\rm d}y}' - f(s)\right)\, {{\rm d}s}
  \\ & \quad 
  + \iot \hh (\ph(s)) \left(\Rcel(s,y') -  \int_{\cal Y} \Rcel(s,y')f(s,y') \,{{\rm d}y}'\right)  f(s)\, {{\rm d}s},
  \label{wf:phenot}
\end{align}
}%
for  every $t \in [0,T],$ and almost everywhere in $ \Omega\times \YY$. 

\end{theorem}

\begin{remark}
We notice that the boundary conditions in \eqref{eq:6} are encompassed within the regularity space $\Hdn$ in \eqref{reg:1:phi} and \eqref{reg:3:sigma}. We also notice that alternative expressions could be considered for some of the terms in the weak formulation of Theorem~\ref{THM:WP}. Specifically, the terms ${\iO} \nabla\ph \cdot \nabla v$ and ${\iO} \nabla\sigma \cdot \nabla v$ could be reformulated using Laplacians taking into account the corresponding zero-flux boundary conditions. Additionally, the duality term involving the time derivative of $\sigma$ may be rewritten as a proper integral, in the light of the time integrability of $\sigma$ as outlined in~\eqref{reg:3:sigma}.
\end{remark}

Under the additional assumption that $\mobm$ and $\hh$ are constant functions, which we set to unity without loss of generality, i.e. 
\begin{align}
\label{ass:additionalmahh}
\mobm \equiv 1, \quad \hh \equiv 1,
\end{align}
uniqueness can also be proved, as established by the following theorem. 
\begin{theorem}
\label{THM:UQ}
Let assumptions \tutte\ and the additional assumptions \eqref{ass:additionalmahh} hold. Then the weak solution $\soluz$ to the initial-boundary value problem  \Sys\ given by Theorem~\ref{THM:WP} is unique. Moreover, if $\{(f_i, \ph_i,\mu_i,\sigma_i)\}_i$, with $i=1,2$, are two solutions to the problem \Sys\ corresponding to the initial data {$\{(f_0^i, \ph_0^i, \sigma_0^i)\}_i$, with $i=1,2$, fulfilling \eqref{ass:data}-\eqref{f0:density:2}}, then the following estimate holds
\begin{align}
	& \non
	\norma{(\ph_1- \ph_2) -(\ph_1- \ph_2)_\Omega}_{\L\infty \Vp}	
	+ \norma{\ph_1- \ph_2}_{\L2 {\Hx1}}
	\\ & \qquad \non
	+ \norma{f_1-f_2}_{C^{1}([0,T]; L^\infty(\Omega) \times L^1(\YY))}
	+ \norma{\sigma_1 - \sigma_2}_{\L\infty {\Lx2} \cap \L2 {\Hx1}}
	\\ & \quad \leq
	K
	( 
	\norma{(\ph_0^1- \ph_0^2) -(\ph_0^1- \ph_0^2)_\Omega}_{\Vp}
	{+ 	\norma{f_0^1 - f_0^2}_{L^\infty(\Omega) \times L^1(\YY)}}
	+ 	\norma{\sigma_0^1 - \sigma_0^2}
	),
	\label{uq:est}
\end{align}
where $K$ is a positive real constant that depends only on $\Omega$, $T$, and on the structure of the problem~\eqref{eq:1}-\eqref{eq:7}.
\end{theorem}

\subsection{Proof of Theorem~\ref{THM:WP} }
\label{SEC:EX}
In proving Theorem \ref{THM:WP}, we proceed formally, acknowledging that, to obtain the estimates below rigorously, we would need to employ a suitable approximation procedure, such as the Faedo--Galerkin approach. This procedure consists in projecting the system onto finite-dimensional subspaces spanned by the eigenfunctions of the elliptic operator subject to homogeneous Neumann boundary conditions. The corresponding orthogonal projections yield a sequence of approximate solutions, and the density of these subspaces in the relevant function spaces ultimately ensures the passage to the limit. However, this is a purely technical step that could be carried out through methods similar to those used, e.g.,  in~\cite{EG, GLS, KS2, RSchS, SS} for related systems, and we thus omit it here.

\step
Key properties of the local cell phenotype distribution

We start by considering the integro-differential equation~\eqref{eq:3} for the local cell phenotype distribution $f$. More specifically, we consider the corresponding integrated-in-time version given by~\eqref{wf:phenot}. Given the uniform boundedness of the function {$\hh$}, under the assumptions in \eqref{prop:M:1}-\eqref{prop:M:2} on ${\MM}$ and the assumptions in~\eqref{prop:R:1} on $\Rcel$, standard theory of integral equations (see e.g. \cite{integrodiff}) ensures that the solution $t \mapsto f(t)$ to the integral equation \eqref{wf:phenot} is in $C^0([0,T])$.

Moreover, under assumptions~\eqref{f0:density:1} on the initial conditions, we have
\begin{align}	
	\label{mass:cons}
	\int_\YY f(x,t,y) \,{{\rm d}y} = \int_\YY f_0(x,y) \,{{\rm d}y}=1 
	\quad 
	\forall t \in [0,T], \, \text{ for a.e. } x \in\Omega,
\end{align}
as it can be easily seen by integrating \eqref{wf:phenot} over the phenotype domain $\YY$. In fact, integrating over $\YY$ yields
\begin{align*}
	& \int_\YY f(x,t,y)\,{{\rm d}y} = \int_\YY f_0 (x,y)\,{{\rm d}y}
	\\ & \quad 
	+ {\theta} \underbrace{\int_\YY \left(\iot \hh (\ph(x,s)) \left( \int_{\cal Y} {\MM}(x,y', y)f(x, s, y') \,{{\rm d}y}' - f(x,s,y)\right)\,{{\rm d}s}\right)\,{{\rm d}y}}_{=0}
	\\ & \quad 
	+ {\underbrace{\int_\YY \left(\iot \hh (\ph(x,s)) \left( \Rcel(x,s,y)- \int_{\cal Y} \Rcel(x,s, y')f(x, s, y') \,{{\rm d}y}' \right)f(x,s,y)\,{{\rm d}s}  \right)\,{{\rm d}y}}_{=0}}=1,
\end{align*}
whence the mass conservation property \eqref{mass:cons} holds. This also entails that 
\begin{align*}
	f(t) \in L^\infty (\Omega) \times L^1(\YY),
	\quad 
	\forall t \in [0,T].
\end{align*}
In view of~\eqref{mass:cons}, to prove that $f$ is a probability density function over the phenotype domain $\YY$, it remains to prove that $f$ is non-negative a.e. on $\Omega \times \YY$. Since this can be proven as similarly done in~\cite{Desv}, here we just summarise the key elements of the method of proof, which relies on a time discretisation argument. Let $n \in \enne$ and, for compactness of notation, omit the dependence on $(x,y)$. Then, we consider the time-discretised problem given by
\begin{align*}
	f^0 &:=f_0,
	\quad	
	\dt f^{n+1} 
	= \Phi(\cdot,f^n)
	-\hh(\cdot) {\Big(\theta - \Rcel + \int_\YY \Rcel f^n \, {\rm d}y \Big)}f^n,
	\\ 
	f^{n+1}(0)& =f_0,
	\quad 
	n \in \enne,
\end{align*}
with 
$\displaystyle{\Phi(r,f) :={\hh(r) {\theta}}\left( \int_\YY {\MM} f \,{{\rm d}y}\right)}$ for every $ r \in \erre$ and $f : \YY \to \erre$, where we recall that, under assumptions~\ref{ass:4:troncante}, the function $\hh(r)$ is non-negative and uniformly bounded for every $r \in \erre$. From this, by employing an inductive method and manipulating the semi-explicit solution directly, one then demonstrates that $f^n$ is non-negative for all $n\in\enne$ and for a.e. $(x,y) \in \Omega \times \YY$. Subsequently, it is a standard matter to prove that, as $n$ tends to infinity, a solution $f$ to the original integro-differential equation is recovered.

A consequence of these properties of $f$, which will be invoked and repeatedly used later on in the proof, is that, for $g \in \{\pp,\qq,\kk\}$,
\begin{align*}
	\Big|\int_\YY g(y) f(t,x ,y) \,{{\rm d}y} \Big|
	\leq 
	\norma{g}_{L^\infty(\YY)}\int_\YY \big|f(t,x, y)\big| \,{{\rm d}y} 
	= \norma{g}_{L^\infty(\YY)} \quad \forall t \in [0,T], \,\, \text{ for  a.e.}\,\, x \in \Omega.
\end{align*}

\step 
Key properties of the local nutrient concentration

We then turn to the initial-boundary value problem \eqref{eq:4}-\eqref{eq:6} for the local nutrient concentration $\sigma$, which can be rewritten as
\begin{align}
\label{eq:revpdesigma}
	\begin{cases}
	\dt \sigma  - \Delta\sigma 
		 = c^* \sigma+b\, \sigma_B 
	\quad &\text{in $Q$,}
	\\ 
	\dn \sigma=0
	\quad &\text{on $\Sigma$,}
	\\ 
	\sigma(0)=\sigma_0
	\quad &\text{in $\Omega$,}
	\end{cases}
\end{align}
with
$\displaystyle{c^*:=  -\hh(\ph) \left(\int_\YY \kk f \,{{\rm d}y}\right) - b}$.
Under assumptions \ref{ass:3:pheno} and \ref{ass:4:troncante}, and given the mass conservation property \eqref{mass:cons}, one can easily show that $\norma{c^*}_{L^\infty(Q)} + \norma{ \sigma_B }_{\L\infty {\Lx2}} \leq C^*$, for a positive and computable constant $C^*$. Hence, through testing with $\sigma$ and then with $-\Delta \sigma$, it is a matter of routine to deduce that
\begin{align*}
	\frac12 \frac {\rm d} {{\rm d}t} \norma{\sigma}_{\Hx1}^2 
	+ \norma{\nabla \sigma}^2
	+ \norma{\Delta \sigma}^2
	\leq 
	C(\norma{\sigma}_{\Hx1}^2+1),
\end{align*}
and then, applying Gr\"onwall's lemma along with elliptic regularity theory, deduce also that
\begin{align*}
	\norma{\sigma}_{\L\infty {\Hx1}
	\cap {\L2 {\Hx2}}}
	\leq C.
\end{align*}
From this, one readily infers, by using equation \eqref{wf:3}, that 
\begin{align*}
	\norma{\dt \sigma}_{\L2 {\Lx2}}
		\leq C (\norma{\Delta\sigma}_{\L2 {\Lx2}} + \norma{\sigma}_{\L2 {\Lx2}})
	\leq C.
\end{align*}

Next, we notice that the \rhs\ of the parabolic equation~\eqref{eq:revpdesigma} is bounded in $\L\infty {\Lx2}$ due to the last assumption in \eqref{ass:data} and the above estimates. Moreover, since we have homogeneous Neumann boundary conditions, and the initial condition $\sigma_0$ is uniformly bounded (cf. the last assumption in \eqref{ass:data}), we invoke the standard parabolic result in \cite[Thm. 7.1,
p. 181]{lady} to conclude that
\begin{align}
	\label{sigma:unif}
	\norma{\sigma}_{L^\infty(Q)}\leq C.
\end{align} 

\step
Energy estimate

Now that we have acquired the above preliminary information, we can proceed with the energy estimate.
We test \eqref{eq:1} with $\mu$, 
\eqref{eq:2} with $-\dt \ph$, 
\eqref{eq:4} with $\sigma$, and then sum together the resulting equalities. Upon rearranging terms, we obtain
\begin{align*}
	& 
	\frac 12\frac {\rm d} {{\rm d}t} \left( 
	\norma{\nabla \ph}^2	
	+2 \iO F(\ph)
	+ \norma{\sigma}^2	
	\right)
	+ \iO \mobm(\ph)|{\nabla \mu}|^2
	+ \norma{\nabla \sigma}^2
	\\ & \quad 
	= 
	\underbrace{\iO \hh(\ph) \left( \sigma \int_\YY \pp f \,{{\rm d}y} 
	-\int_\YY \qq f \,{{\rm d}y} \right) \mu}_{=:I_1}
	\\ & \qquad 
	\underbrace{-\iO \hh(\ph)\sigma \left(\int_\YY \kk f \,{{\rm d}y} \right)\sigma
  + b \iO \bigl(\sigma_B - \sigma\bigr) \sigma}_{=:I_2}.
\end{align*}
Now, to control $I_2$ we resort to \eqref{mass:cons}, \Holder's and Cauchy--Schwarz inequalities, the boundedness of $\hh$ and $\kk$, and the regularity of $\sigma_B$.
Simple computations, using also the uniform bound for $\sigma$ in \eqref{sigma:unif}, lead to 
\begin{align*}
	|I_2| 
	& \leq
	\norma{\hh}_{L^\infty(\Omega)}\norma{\sigma}_{L^\infty(\Omega)}
	\iO \norma{\kk}_{L^\infty(\YY)}\Big|\int_\YY f \,{{\rm d}y} \Big|
	+ C(\norma{\sigma_B}^2 +1)
	\leq C	.
\end{align*}
To bound $I_1$, we need to carry out a few steps. First, we notice that 
\begin{align*}
	|I_1| = \Big|\iO \hh(\ph) \left( \sigma \int_\YY \pp f \,{{\rm d}y} 
	-\int_\YY \qq f \,{{\rm d}y} \right) \mu \anold{\Big|}
	\leq C \iO |g | | \mu|,
\end{align*} 
where $\displaystyle{g:= \hh(\ph) \left( \sigma \int_\YY \pp f \,{{\rm d}y} 
	-\int_\YY \qq f \,{{\rm d}y} \right)}$ is uniformly bounded in view of assumptions \ref{ass:3:pheno} and \ref{ass:4:troncante}, and the bound for $\sigma$ in \eqref{sigma:unif}.
In more detail, using also the Poincar\'e inequality, we find
\begin{align*}
	|I_1|
	\leq \iO|g|| \mu-\mu_\Omega| + |\Omega| |g| |\mu_\Omega|
	\leq 
	\frac {m^*}2 \norma{\nabla \mu}^2 + C + C |\mu_\Omega|.
\end{align*}
Combining the above estimates, we obtain 
\begin{align}
	& 
	\frac 12\frac {\rm d} {{\rm d}t} \left( 
	\norma{\nabla \ph}^2	
	+ 2 \iO F_1(\ph)
	+ \norma{\sigma}^2	
	\right)
	+ \frac {m^*}2\norma{\nabla \mu}^2
	+ \norma{\nabla \sigma}^2
	\leq 
	C + C |\mu_\Omega|.
		\label{energy:parz}
\end{align}

We are left to bound the last term on the \rhs\ of \eqref{energy:parz}. To this end, we test equation \eqref{eq:2} with ${|\Omega|}^{-1}$, and use the assumptions in \eqref{growth:F:1}-\eqref{growth:F:2} to infer that
\begin{align}	
	\label{mean}
	|\mu_\Omega|
	\leq C \iO (|F'_1(\ph)| + 1+ |\ph|^2)
	= 
	C (1+ \norma{F_1(\ph)}_1+ \norma{\ph}^2).
\end{align}
We now sum together \eqref{energy:parz} and \eqref{mean}, and then use assumption \ref{ass:1:potenziale} to find 
\begin{align*}
\non
	& 
	\frac 12\frac {\rm d} {{\rm d}t} \left( 
	\norma{\ph}^2_V	
	+ 2 \norma{F_1(\ph)}_1
	+ \norma{\sigma}^2	
	\right)
	+ \frac {m^*}2\norma{\nabla \mu}^2
	+ \norma{\nabla \sigma}^2
	\leq C(
\norma{F_1(\ph)}_1+ \norma{\ph}^2 +1),
\end{align*}
which, by applying Gr\"onwall's lemma, leads to
\begin{align*}
	& \norma{\ph}_{\L\infty {\Hx1} }
	+ \norma{F_1(\ph)}_{\L\infty {\Lx1} }
	+ \norma{\nabla \mu}_{\L2 {\Lx2}}
	\leq C.
\end{align*}
Returning to equation \eqref{mean}, one sees that $\norma{\mu_\Omega}_{L^2(0,T)} \leq C$ and thus, using Poincaré inequality, one infers that
\begin{align*}
	\norma{F'(\ph)}_{\L2 {\Lx2}}
	+\norma{\mu}_{\L2 {\Hx1}}
	\leq C.
\end{align*}
Moreover, using the weak formulation of equation \eqref{wf:1}, we find
\begin{align*}
	\norma{\dt \ph}_{\L2 \Vp}
	\leq C.
\end{align*}
Finally, using the aforementioned bounds and elliptic regularity theory, we infer from equation \eqref{wf:1} that
\begin{align*}
	\norma{\ph}_{\L2 {\Hx2}}
	\leq C.
\end{align*}
This concludes the proof of Theorem~\ref{THM:WP}. \qed

\subsection{Proof of Theorem~\ref{THM:UQ} }
\label{SEC:UQ}
The method for proving uniqueness relies on a conventional strategy that involves demonstrating the stability estimate outlined in \eqref{uq:est}, from which it follows that uniqueness is ensured when considering solutions generated from the same initial data.

To begin with, we introduce the following abridged notation for the differences between two solutions and the corresponding initial data:
\begin{align*}
	\ph & := \ph_1-\ph_2,
	\quad 
	\mu:=\mu_1-\mu_2,
	\quad 
	f  := f_1- f_2
	\quad  
	\sigma :=\sigma_1-\sigma_2,
	\\
	\ph_0 & := \ph_0^1-\ph_0^2,
	\quad  
	{f_0 :=f_0^1-f_0^2},
	\quad  
	\sigma_0 :=\sigma_0^1-\sigma_0^2.
\end{align*}
Then, we notice that the above differences $(f,\ph,\mu,\sigma)$ solve the following system, which we write in strong form for convenience:
\begin{alignat}{2}
	& \dt \ph -  \Delta \mu
	\label{cd:eq:1}
	= \Rold(f, f_1, \sigma)
	\quad && \text{in $Q$,}
	\\
	\label{cd:eq:2}
	& \mu = - \Delta \ph + (F' (\ph_1)-F' (\ph_2))
	\qquad && \text{in $Q$,}
	\\
		& {\dt f = \theta \, \left(\int_{\cal Y} {\MM} f \,{{\rm d}y} - f\right)}
		\label{cd:eq:3}
	 - {\left( \int_\YY \Rcel  f\,{\rm d}y\right) f_1 +\left( \int_\YY \Rcel  f_2 \,{\rm d}y\right) f}
	\qquad && \text{in $Q \times \YY$,} 
	\\ 
	\label{cd:eq:4}
	& \dt \sigma - \Delta \sigma
	+ \sigma \left(\int_\YY \kk  f_1  \,{{\rm d}y}\right)
	+ \sigma_2 \left(\int_\YY \kk  f \,{{\rm d}y}\right)	
	=
	- b \sigma
	\quad && \text{in $Q$,}
	\\
	\label{cd:eq:5}
	& \dn \ph= \dn \mu = \dn \sigma
	=
	0
	\quad && \text{on $\Sigma$,}
	\\
	\label{cd:eq:6}
	& \ph(0)=\ph_0, \quad \sigma(0)=\sigma_0
	\quad && \text{in $\Omega$,} \phantom{\int_\YY}
	\\
	\label{cd:eq:7}
	& f(0)={f_0}
	\quad && \text{in $\Omega \times \YY$,}	
\end{alignat}
where we have
\begin{align*}
	\Rold(f, f_1, \sigma):=\sigma \left(\int_\YY \pp  f_1  \,{{\rm d}y} \right)
	+ \sigma_2 \left(\int_\YY \pp  f  \,{{\rm d}y} \right)
	-\int_\YY \qq  f  \,{{\rm d}y}.
\end{align*}
Equation \eqref{cd:eq:3} is linear and, since we have chosen $\hh\equiv 1$, we can then adopt the same strategy as the one employed in the proof of Theorem~\ref{THM:WP}, which makes it possible to obtain results on $f:=f_1-f_2$ that are analogous to those obtained on the local cell phenotype distribution. Namely, we have that there exists a unique solution to \eqref{cd:eq:3} such that:
\begin{align*}
	& f \in C^0 ([0,T]; L^\infty (\Omega) \times L^1(\YY)),
	\\ &\text{and}\quad
	\exists \,C>0 :  \quad 
	\norma{f_1-f_2}_{C^0 ([0,T]; L^\infty (\Omega) \times L^1(\YY))}\leq C {\norma{f_0^1-f_0^2}_{L^\infty(\Omega) \times L^1(\YY)}}.
\end{align*}

Next, testing \eqref{cd:eq:1} with $\dfrac 1{|\Omega|}$, we find
\begin{align}
 	\non
 	\frac {\rm d} {{\rm d}t} \ph_\Omega 
 	& =
 	(\Rold(f, f_1, \sigma))_\Omega
 	\\ &  \non
 	= \frac 1{|\Omega|}\iO  \sigma \left(\int_\YY \pp f_1 \,{{\rm d}y} \right)
	+ \frac 1{|\Omega|}\iO \sigma_2 \left(\int_\YY \pp f \,{{\rm d}y} \right)
 	\\ &  \quad 
	- \frac 1{|\Omega|}\iO \int_\YY \qq f \,{{\rm d}y}.
	\label{reh:mean}
\end{align}
Recalling the properties of the operator $\cal N$ summarised in Section \ref{SEC:NOT}, we first subtract the above from \eqref{cd:eq:1} and test the resulting equation with ${\cal N}(\ph - \ph_\Omega)$, next we test \eqref{cd:eq:2} with $-(\ph- \ph_\Omega)$, then test \eqref{cd:eq:5} with $\sigma$, and finally sum together the resulting equations to obtain
\begin{align}
	\non & \frac 12 \frac {\rm d} {{\rm d}t} \norma{\ph-\ph_\Omega}^2_{*}
	+ \norma{\nabla \ph}^2
	+ \iO (F' (\ph_1)-F' (\ph_2)) (\ph- \ph_\Omega)
	+\frac 12 \frac {\rm d} {{\rm d}t} \norma{\sigma}^2
	+ \norma{\nabla \sigma}^2
	\\ & \quad  \non
	= \iO \big(\Rold(f, f_1, \sigma) - (\Rold(f, f_1, \sigma))_\Omega\big) {\cal N}(\ph-\ph_\Omega)
	- b \norma{\sigma}^2
	\\ & \qquad \label{cd:est}
	- \iO \sigma \left(\int_\YY \kk f_1 \,{{\rm d}y}\right)\sigma
	- \iO \sigma_2 \left(\int_\YY \kk f \,{{\rm d}y}\right)	\sigma.
\end{align}
Note that, due to assumptions \ref{ass:1:potenziale} and classical interpolation estimates, it holds that
\begin{align*}
	\iO (F' (\ph_1)-F' (\ph_2)) (\ph- \ph_\Omega)\geq - \frac 14 \norma{\nabla \ph}^2 - C \norma{\ph-\ph_\Omega}^2_{*},
\end{align*}
where the first term on the \rhs\ can be absorbed, whereas the other one can be moved to the left-hand side. Moreover, using Young's inequality, along with some interpolation and the regularity of  $f_1$, we infer that
\begin{align*}
	& \iO \big(\Rold(f, f_1, \sigma) - (\Rold(f, f_1, \sigma))_\Omega\big) {\cal N}(\ph-\ph_\Omega)
	- b \norma{\sigma}^2
	\\ & \quad 
	\leq 
	\norma{\Rold(f, f_1, \sigma) - (\Rold(f, f_1, \sigma))_\Omega}_{*}\norma{\ph-\ph_\Omega}
	+C \norma{\sigma}^2
	\\ &  \quad 
	\leq 
	\frac 14 \norma{\nabla \ph}^2
	+ C \norma{\ph-\ph_\Omega}^2_{*}
	+ C \norma{\sigma}^2 .
\end{align*}
Finally, using Young's inequality along with the regularity of $f_1$ and $\sigma_2$, the last line of \eqref{cd:est} can be controlled as follows:
\begin{align*}
	& 
	- \iO \sigma \left(\int_\YY \kk f_1 \,{{\rm d}y}\right)\sigma
	- \iO \sigma_2 \left(\int_\YY \kk f \,{{\rm d}y}\right)	\sigma
	\leq C \norma{\sigma}^2,
\end{align*}
so that, employing Gr\"onwall's lemma, we obtain \eqref{uq:est}. This concludes the proof of Theorem~\ref{THM:UQ}. \qed

\section{Numerical simulations}
\label{SEC:NUM}

\subsection{Numerical scheme}
\label{SUBSEC:methods}

We consider the following finite element space
\begin{equation*}
V_h = \{\chi \in C^0(\ov{\Omega}):\,\,\chi |_{K} \in \mathcal{P}^1(K)\,\,\forall\,K \in \mathcal{T}_h\} \subset H^1(\Omega),
\end{equation*}
where $\mathcal{T}_h$ is a partition of the spatial domain $\Omega$ into triangular elements, and we denote by $\mathcal{P}^1(K_j)$ the space of polynomials of order one on the tetrahedron $K_j \in \mathcal{T}_h$. Moreover, we consider a uniform discretisation of the phenotype domain $\YY := [y^m, y^M]$ of step $\Delta y = (y^M- y^m)/N_y$, with $N_y \in\enne$, whereby $y_i = y^m + i\Delta y$ for $i = 0, \ldots, N_y$. Similarly, we consider a uniform discretisation of the time domain $[0,T]$ of step $\Delta t = T/N_t$, with $N_t \in\enne$, whereby $t^n = n \Delta t$ for $n = 0, \ldots ,N_t$. Finally, we denote the set of initial data by $(\varphi_h^0, \mu_h^0, \sigma_h^0, \textbf{f}_h^{\text{ }0}) \in V_h \times V_h \times V_h \times V_h^{N_y}$.

Building on the method employed in~\cite{pozzi2022t}, assuming the model function $\mobm$ to be constant, i.e. $\mobm(\ph) \equiv m \in \mathbb{R}^+$, we solve numerically the system~\eqref{eq:1}-\eqref{eq:4} through the following scheme. For any $n = 1, \ldots, N_t$, find $(\varphi_h^{n+1},\mu_h^{n+1},\sigma_h^{n+1}, \textbf{f}_h^{\text{ }n+1}) \in V_h \times V_h \times V_h \times V_h^{N_y}$ such that  for all $(v_h, q_h , w_h, \textbf{g}_h)\in V_h \times V_h \times V_h \times V_h^{N_y}$:
\begin{equation}
\Big( \frac{\varphi^{n+1}_h - \varphi^{n}_h}{\Delta t}, v_h \Big) =  - m \, (\nabla \mu^{n+1}_h, \nabla v_h)
	+ (\hh(\ph^{n}_h)(\sigma_h^{n}\ov{\mathcal{P}} - \ov{\mathcal{Q}}), v_h),
	\label{eq::phi_wf}
\end{equation}
\begin{equation}
(\mu^{n+1}_h, \textit{q}_h) = \eps (\nabla \varphi^{n+1}_h, \nabla \textit{q}_h) + \frac{1}{\eps} (\anold{F}_c'(\varphi^{n+1}_h), \textit{q}_h) 
+  \frac{1}{\eps}(\anold{F}_e'(\varphi^{n}_h), \textit{q}_h),
	\label{eq::mu_wf}
\end{equation}
\begin{equation}
	\Big(\frac{\textbf{f}_h^{\text{ }n+1} - \textbf{f}_h^{\text{ }n}}{\Delta t}, \textbf{g}_h \Big) = 
	\alpha \theta (\anold{\hh}(\varphi^{n}_h)\anold{\ov{\bf{M}}}, \textbf{g}_h) - \alpha \theta (\anold{\hh}(\varphi^{n}_h)\textbf{f}_h^{\text{ }n}, \textbf{g}_h)
	+\alpha (\anold{\hh}(\varphi^{n}_h)(\mathcal{R} - \ov{\mathcal{R}})\textbf{f}_h^{\text{ }n}, \textbf{g}_h),
	\label{eq::f_wf}
\end{equation}
\begin{equation}
	\Big(\frac{\sigma^{n+1}_h - \sigma^{n}_h}{\Delta t}, \textit{w}_h \Big) = - D_{\sigma} (\nabla \sigma^{n+1}_h, \nabla w_h)
	- (\anold{\hh}(\varphi^{n}_h)\sigma_h^{n+1}\ov{\mathcal{K}}, w_h) + (b(\sigma_B - \sigma_h^{n+1}), w_h) .
	\label{eq::sigma_wf}
\end{equation}
In~\eqref{eq::phi_wf}-\eqref{eq::sigma_wf}, the notation ($\cdot$, $\cdot$) is used for the standard $L^2$-inner product over $\Omega$, and the following splitting of the Cahn--Hilliard potential is prescribed
\begin{equation*}
	F_c(\varphi^{n+1}_h) := \frac{(\varphi^{n+1}_h)^4 +1}{4}, \qquad F_e(\varphi^{n}_h) := -\frac{(\varphi^{n}_h)^2}{2},
\end{equation*}
so as to ensure that the numerical scheme is gradient stable~\cite{tierra2015numerical}. Moreover, 
$$
\hh(\varphi^n_h) := (1 + \varphi_h^n)/2,
$$
and $\ov{\mathcal{P}} \in V_h$, $\ov{\mathcal{Q}} \in V_h$, and $\ov{\mathcal{K}} \in V_h$ are the numerical approximations of the integrals 
$$
\int_\YY \pp(y) f(x,t,y) \,{{\rm d}y}, \qquad \int_\YY \qq(y) f(x,t,y) \,{{\rm d}y}, \qquad \int_\YY \kk(y) f(x,t,y) \,{{\rm d}y}, 
$$
while $\ov{\mathcal{R}} \in V_h$ and $\ov{\bf{M}} \in V_h^{N_y}$ are the numerical approximations of the integrals 
$$
\int_\YY \Rcel(x,t,y) f(x,t,y) \,{{\rm d}y}, \qquad \int_{\cal Y} {\MM}(x,y',y)f(x,t,y') \,{{\rm d}y}' ,
$$
which are all computed through the corresponding Riemann sum with partition $y_0, \ldots, y_{N_y}$.

\subsection{Set-up of numerical simulations}
\label{SUBSEC:setup}
To carry out numerical simulations, we let the spatial domain $\Omega$ be a square of side $1$, which we then partition into $2,88 \times 10^4$ uniformly distributed triangular elements. Moreover, we choose $y^m = 0$, $y^M = 2$, $N_y = 17$, and $\Delta t = 10^{-3}$, and perform numerical simulations using the open-source computing platform $\mathtt{FEniCS}$~\cite{alnaes2015fenics, logg2012dolfin}. We preliminary verified that, under this setting, the numerical scheme employed ensures that $\textbf{f}_h^{\text{ }n}$ satisfies properties~\eqref{proprel:f}.

As for the model parameters, we set $\eps = 10^{-2}$, $D_{\sigma} = 10^2$, and $b = 10^4$, while we tune the values of the parameters $\alpha$ and $\theta$ as detailed in the next subsection. Moreover, we set $\sigma_B \equiv s_B = 1$ and we use the following definitions of the model functions
$$
\mobm(\ph) \equiv m = 10^{-2}, \qquad \pp(y) \equiv p = 1.5, \qquad \qq(y) := (1 - y)^2 , \qquad \kk(y) \equiv k = 1,
$$
$$
\mathcal{R}(x,t,y) \equiv \mathcal{R}(y) := 1 - 0.1(1-y)^2, \quad \mathcal{M}(x,y',y) \equiv \mathcal{M}(y',y) := \dfrac{\exp[-100(y' - y)^2]}{\displaystyle{\int_{y^m}^{y^M} \exp[-100(y' - y)^2] {\rm d}y}}.
$$
Note that, under these definitions, the phenotype $y=1$ is the maximum point of the fitness function $\mathcal{R}$ for all $(x,t)$ (i.e. the fittest phenotype), which is also the maximum point of the net cell  proliferation rate $\sigma(x,t) \, \pp(y) - \qq(y)$ for all $(x,t)$. Note also that the kernel $\mathcal{M}(y',y)$ satisfies the normalisation condition $\int_{y^m}^{y^M} \mathcal{M}(\cdot,y) {\rm d}y = 1$.

As for the initial conditions, assuming that at $t=0$ the tumour occupies the region $\Omega^T_0 := \{x \in \Omega: (x - x_c)^2 \leq 0.01  \}$, with $x_c = (0.5, 0.5)$ being the centre of the tumour, while the remaining part of the spatial domain is initially occupied by healthy cells only, we set 
$$
\varphi_0:=
\begin{cases}
1, \qquad \text{on } \; \Omega^T_0,
\\
-1, \,\quad \text{on } \; \Omega \setminus \Omega^T_0.
\end{cases}
$$
Furthermore, assuming the majority of cells to have initially phenotype $y=\bar{y}_0$ throughout the spatial domain $\Omega$, we set 
\begin{equation}
\label{def:f0}
f_0(x,y) \equiv f_0(y) := \dfrac{\exp[-a \, (y - \bar{y}_0)^2]}{\displaystyle{\int_{y^m}^{y^M} \exp[- a (y - \bar{y}_0)^2] {\rm d}y}},
\end{equation}
and we tune the values of the parameters $a \in \mathbb{R}^+$ and $\bar{y}_0 \in [y^m,y^M]$ as detailed in the next subsection. Note that, under definition~\eqref{def:f0}, the distribution $f_0$ satisfies the normalisation condition $\int_{y^m}^{y^M} f_0(\cdot,y) {\rm d}y = 1$.

Finally, given the above definitions of $\varphi_0$ and $f_0$, we define $\sigma_0$ as the corresponding solution to the steady-state version of~\eqref{eq:4}, that is, the elliptic equation
$$
- D_{\sigma} \, \Delta \sigma_0 + \hh(\ph_0) \, \sigma \, \left(\int_{\mathcal{Y}} \kk f_0 \,{{\rm d}y}\right) = b\left(\sigma_B - \sigma\right) \qquad \text{in $Q$},
$$
subject to zero-flux boundary conditions, which we compute numerically.

\subsection{Main results of numerical simulations}
\label{SUBSEC:results}
We first carry out numerical simulations choosing $\alpha = 5 \times 10^{2}$, considering different values of the rate of phenotypic changes $\theta$, i.e. $\theta \in \{0, 0.3, 0.5, 0.7\}$, and defining the initial phenotype distribution $f_0$ via~\eqref{def:f0} with $a=2.5$ and $\bar{y}_0=1.75$. The results obtained, which are summarised in Figure~\ref{fig::1}, show that the tumour, i.e. the region
\begin{equation}
\label{def:OmegaT}
\Omega^{\rm T}(t) :=  \left\{x \in \Omega : \hh(\ph(x,t)) > 0 \right\},
\end{equation}
expands over time preserving radial symmetry (cf. the bottom panels of Figure~\ref{fig::1}). Moreover, as indicated by the results of numerical simulations displayed in the top panels of Figure~\ref{fig::1}:
\begin{itemize}
\item[(i)]  For all $t>0$, the local cell phenotype distribution $f$ remains of a Gaussian-like form throughout $\Omega$.
\item[(ii)] As $t \to \infty$, the local mean phenotype, defined via~\eqref{def:meanvarphen}, converges to the fittest phenotype $y=1$ throughout $\Omega^{\rm T}(t)$.
\item[(iii)] As $t \to \infty$, there is a proportional relationship between the local phenotype variance, defined via~\eqref{def:meanvarphen}, and the rate of phenotypic changes $\theta$, throughout $\Omega^{\rm T}(t)$ -- i.e. as $\theta$ gets smaller, $f(x,t,y)$ becomes concentrated as an increasingly sharp Gaussian-like distribution centred at $y=1$ for all $x \in \Omega^{\rm T}(t)$ as $t \to \infty$.
\end{itemize}
Biologically, these numerical results support the ideas that the rate at which tumour cells undergo phenotypic changes, $\theta$, impacts on the level of intra-tumour phenotypic heterogeneity, and the prevailing phenotype in the tumour will ultimately be the fittest phenotype. This is in line with the findings of earlier theoretical studies on the evolutionary dynamics of cancer cells in solid tumours~\cite{Chisholmetal2016,Marusyketal2012,Michor2010,villaetal2021b}. We verified that, as it could be expected, tuning the value of $\alpha$ impacts only on the time scale over which the dynamics of the cell phenotype distribution $f$ take place, while leaving the key qualitative features (i)-(iii) intact.

\begin{figure}
    \centering
	\includegraphics[width=0.32\textwidth, trim=30 5 35 15,clip]{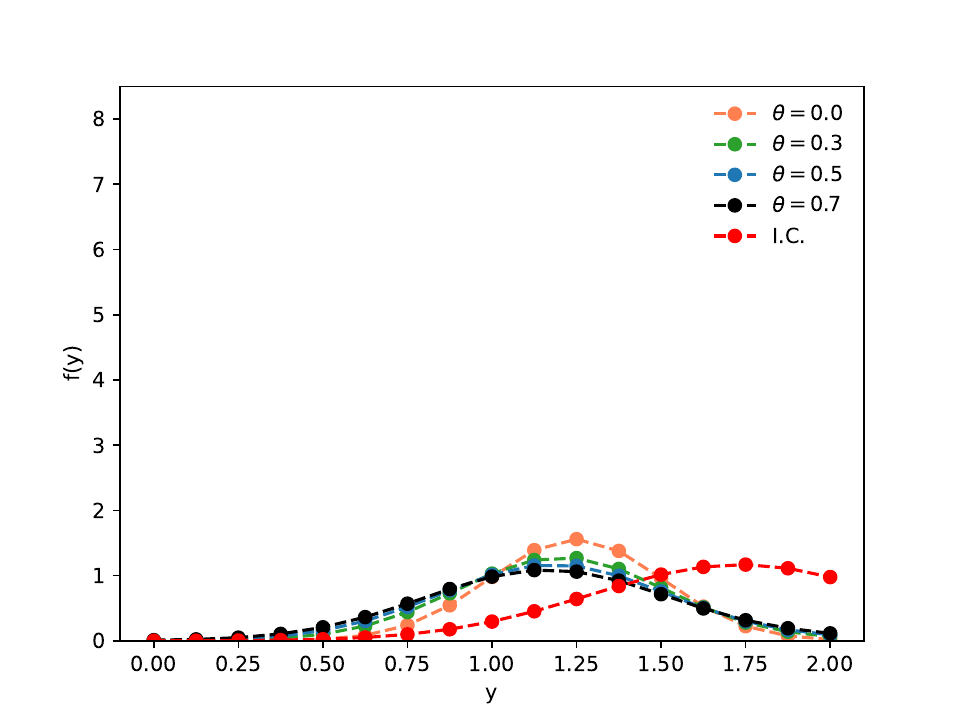}
	\includegraphics[width=0.32\textwidth, trim=30 5 35 15,clip]{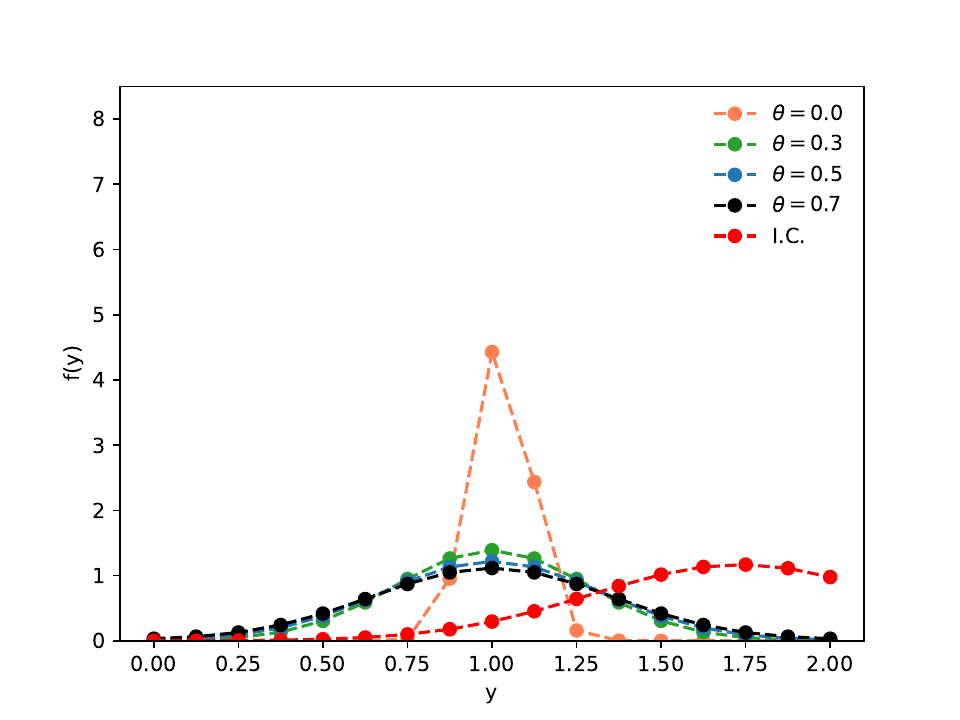}
	\includegraphics[width=0.32\textwidth, trim=30 5 35 15,clip]{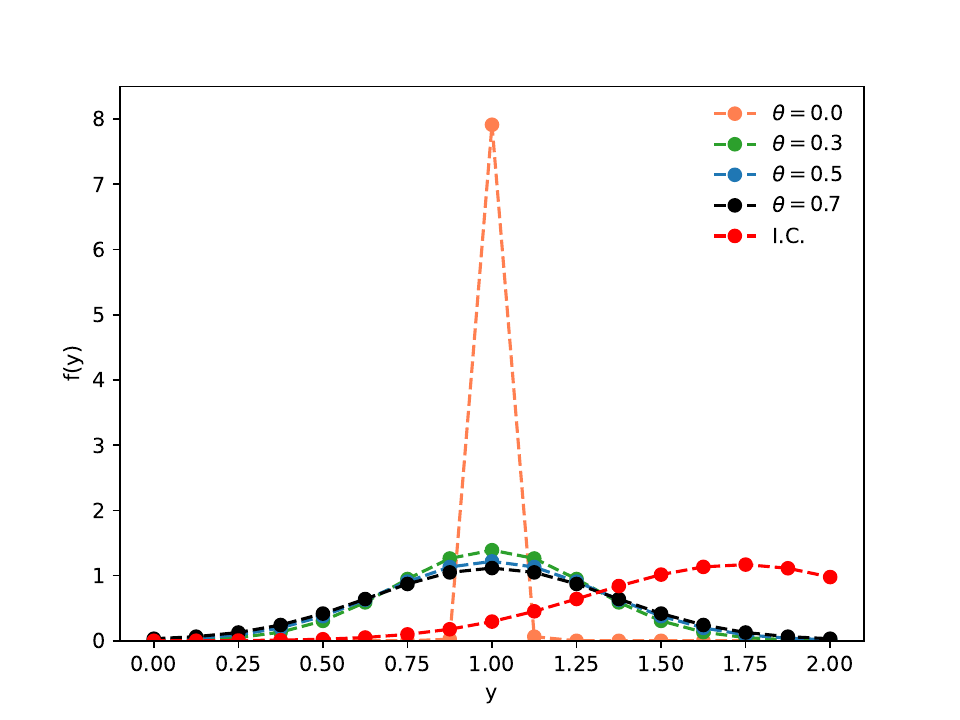}
	\includegraphics[width=0.3\textwidth]{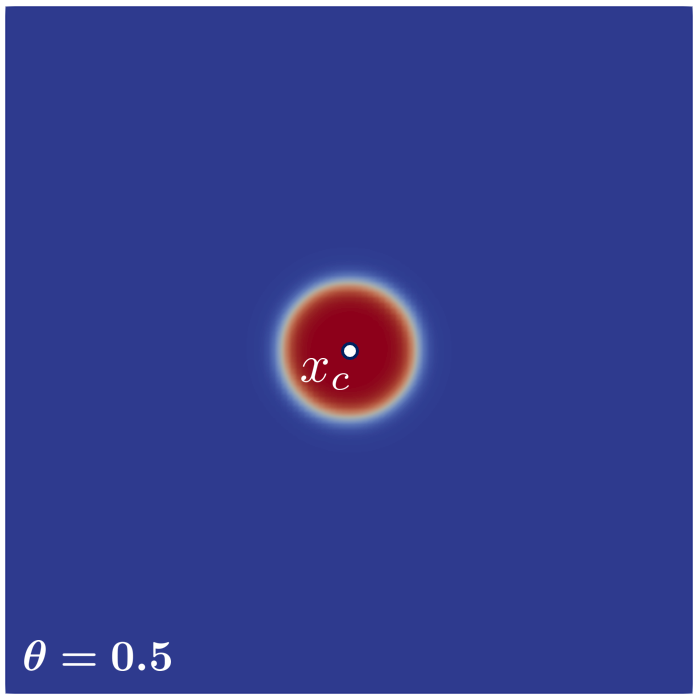}\hfil
	\includegraphics[width=0.3\textwidth]{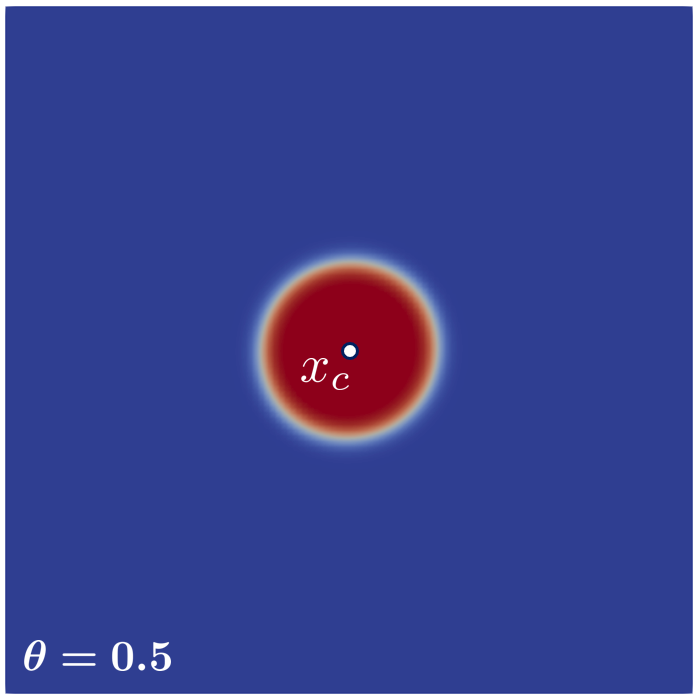}\hfil
	\includegraphics[width=0.3\textwidth]{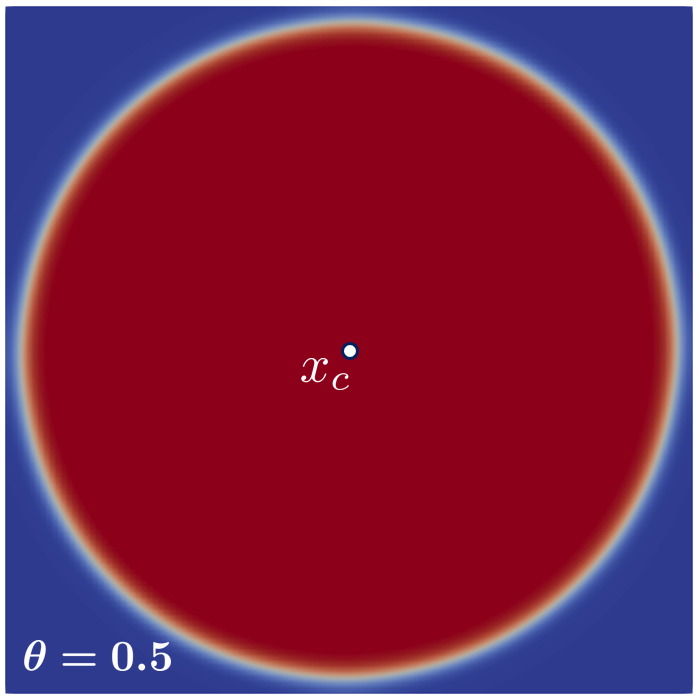}\\
	\includegraphics[width=0.27\textwidth]{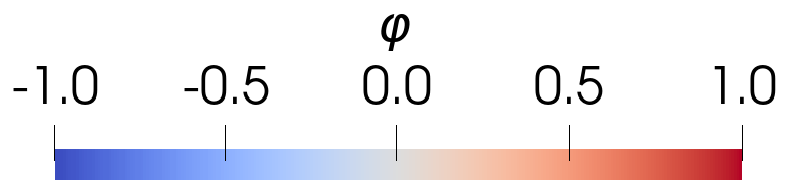}
    \caption{ \small {\bf Top panels.} Plots of the local cell phenotype distribution $f$ at $x=(0.5, 0.5)$ (i.e. the centre of the tumour $x_c$ displayed in the bottom panels), at $t = 0.1$ (left), $t = 1.3$ (centre), and $t = 5.4$ (right), for different values of the model parameter $\theta$ specified in the legend. The initial local cell phenotype distribution $f_0$ is plotted in red for reference (cf. I.C. in the legend). {\bf Bottom panels.} Plots of the phase variable $\varphi$ at $t = 0.1$ (left), $t = 1.3$ (centre), and $t = 5.4$ (right) for $\theta = 0.5$. Qualitatively similar dynamics of $\varphi$ are obtained for the other values of $\theta$ considered here (results not shown).}
    \label{fig::1}
\end{figure}

A corroboration of (i) and (ii) is shown in Figure~\ref{fig::3} for the case when $\theta=0.5$ -- qualitatively similar results are obtained for the other values of $\theta$ considered here (results not shown). In Figure~\ref{fig::3}, the numerical solution of~\eqref{eq:3} is seen to remain of a Gaussian-like form for all $(x,t)$, with the local mean phenotype that starts moving from $\bar{y}_0=1.75$ towards $y=1$ as soon as $\varphi(x,t)$ grows above $-1$ (i.e. as soon as some tumour cells become present at position $x$), so that $\hh(\ph(x,t))$ grows above $0$ making the left-hand side of equation~\eqref{eq:3} to become non-zero at $(x,t)$, thus initiating the evolution of the local cell phenotype distribution $f(x,t,y)$. 
\begin{figure}
    \centering
    \includegraphics[width=0.3\textwidth, trim=25 5 45 15,clip]{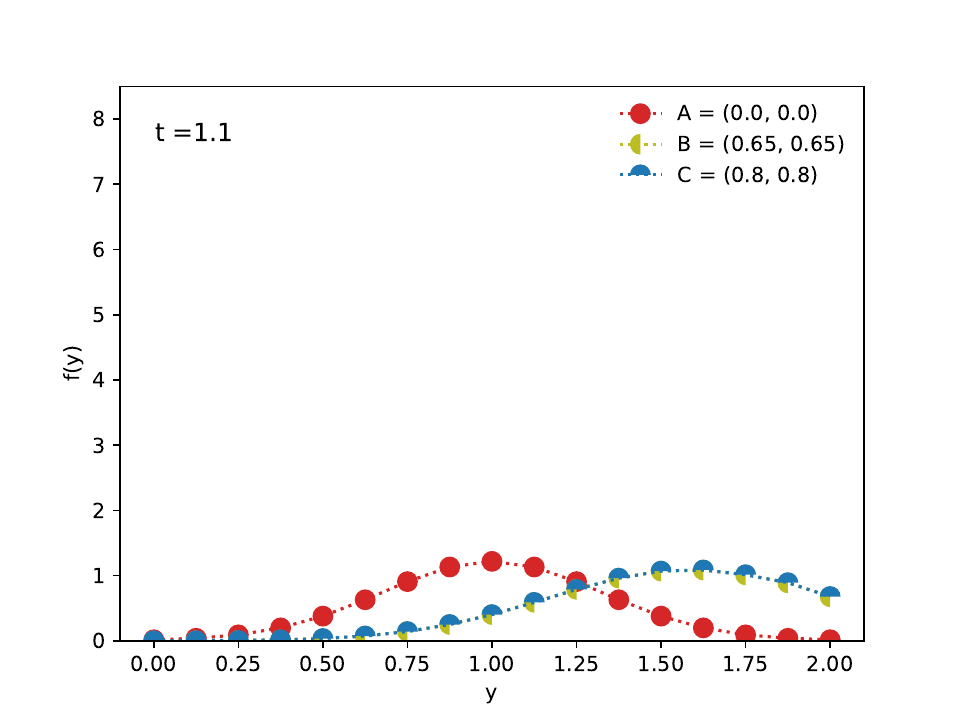}
	\includegraphics[width=0.3\textwidth, trim=25 5 45 15,clip]{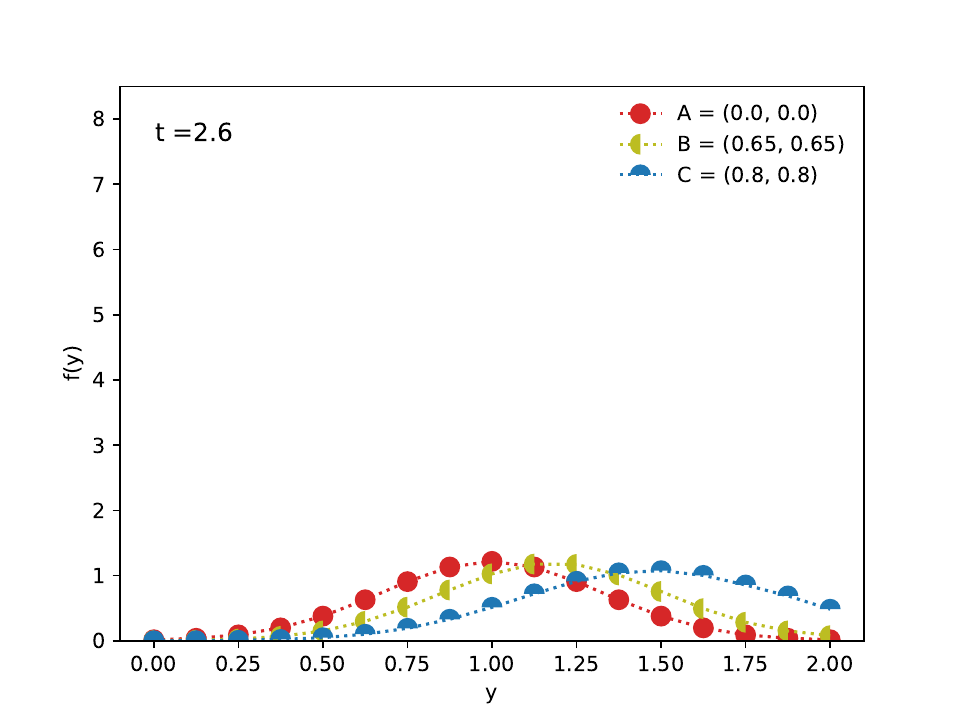}
	\includegraphics[width=0.3\textwidth, trim=25 5 45 15,clip]{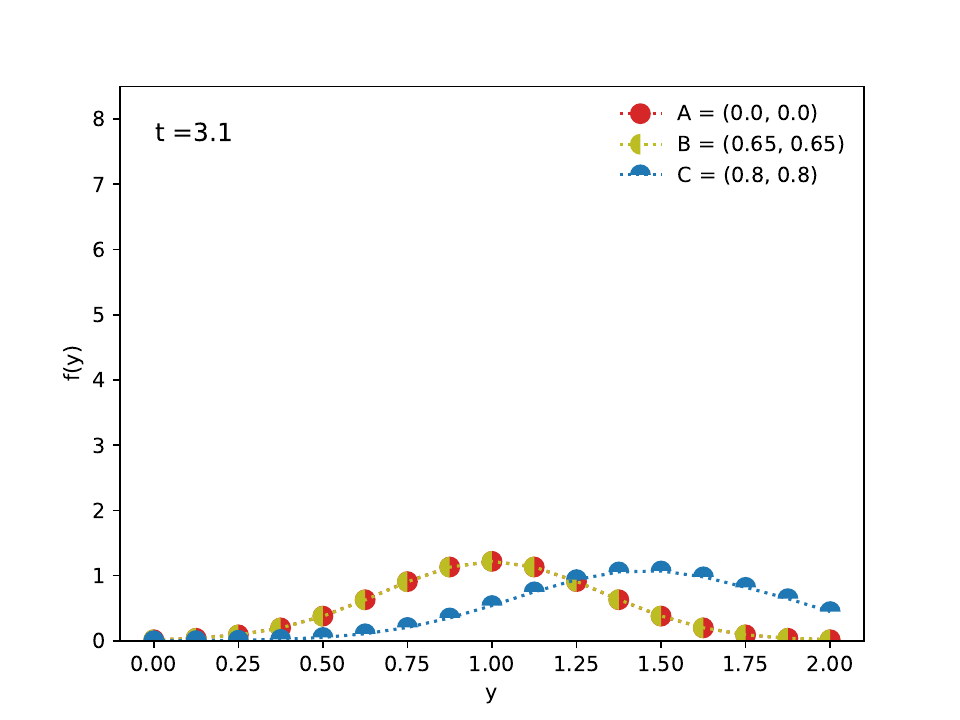}
	\includegraphics[width=0.3\textwidth, trim=25 5 45 15,clip]{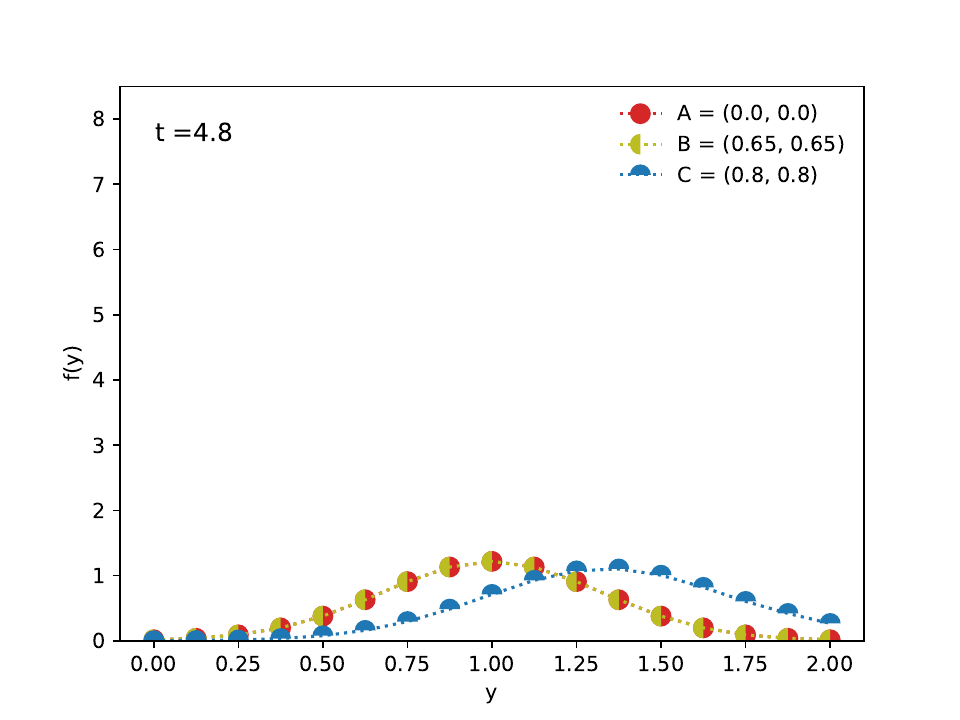}
	\includegraphics[width=0.3\textwidth, trim=25 5 45 15,clip]{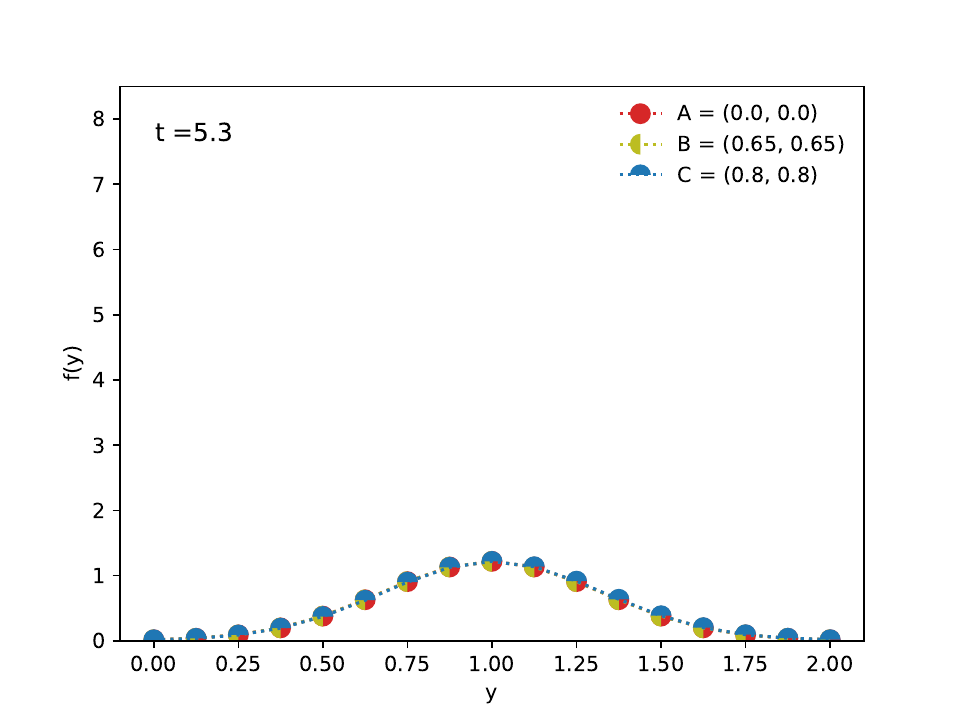}
	\includegraphics[width=0.231\textwidth, trim=0 -10 0 0,clip]{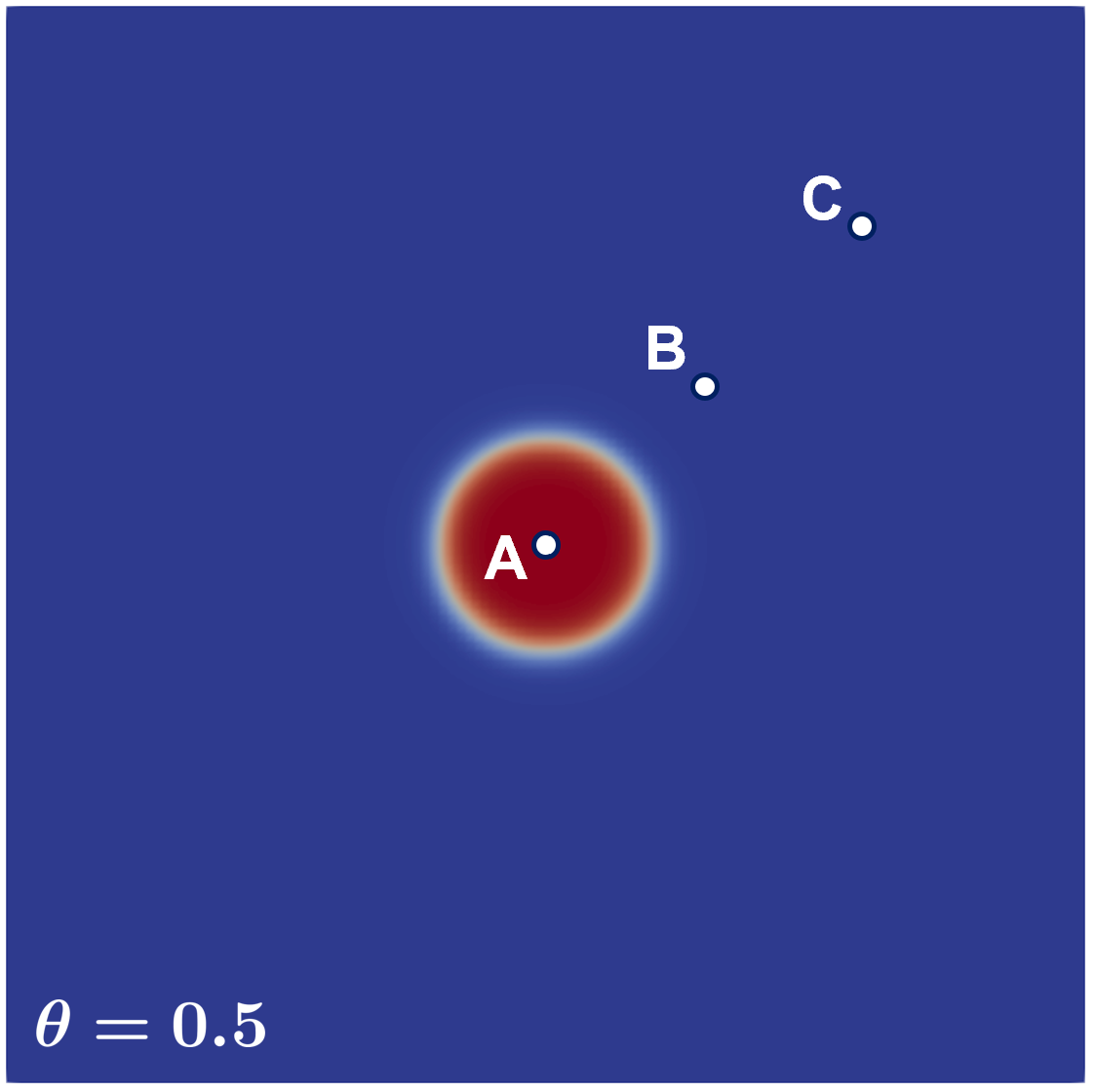}
	\includegraphics[width=0.05\textwidth]{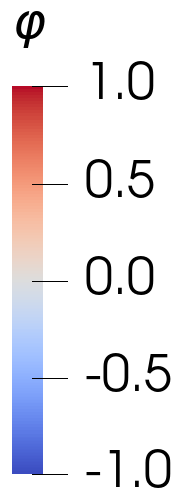}
    \caption{\small {\bf Top left-bottom centre panels.} Plots of the local cell phenotype distribution $f$ at the points $x = (0.5, 0.5)$ (red lines), $x= (0.65, 0.65)$ (green lines), and $x= (0.8, 0.8)$ (cyan lines) (i.e. the points ${\rm A}$, ${\rm B}$, and ${\rm C}$ displayed in the bottom right panel), at $t=1.1$ (top left), $t=2.6$ (top centre), $t=3.1$ (top right), $t=4.8$ (bottom left), and $t=5.3$ (bottom centre), for $\theta = 0.5$. {\bf Bottom right panel.} Plot of the initial phase variable $\varphi_0$.}
    \label{fig::3}
\end{figure}

In order to investigate, in controlled scenarios, the impact of the shape of the local cell phenotype distribution $f$ on the expansion rate of the region $\Omega^{\rm T}(t)$ defined via~\eqref{def:OmegaT} (i.e. to investigate how the phenotypic composition of the tumour affects tumour growth), we also carry out numerical simulations setting $\alpha=0$, so that $f(x,t,y) = f_0(x,y)$ for all $t>0$, defining $f_0$ via~\eqref{def:f0} with $a= 2.5$ (i.e. $f_0$ is defined as a sharp Gaussian-like distribution centred at $\bar{y}_0$) and $\bar{y}_0=1$ (to which we refer as initial condition IC2) or $\bar{y}_0=1.7$ (to which we refer as initial condition IC1). We then compare the results obtained in such controlled scenarios with those previously obtained in the case where $\alpha > 0$ (i.e. $\alpha = 5 \times 10^{2}$), for $\theta=0.5$ and $f_0$ defined via~\eqref{def:f0} with $a=2.5$ and $\bar{y}_0=1.75$ (to which we refer as initial condition IC0). As summarised by the plots in Figure~\ref{fig::2}, the region $\Omega^{\rm T}(t)$ expands more quickly in the controlled scenario where the local cell phenotype distribution remains concentrated around the fittest phenotype $y=1$ (purple curves, corresponding to the initial condition IC2 in Figure~\ref{fig::2}), as this correlates with a relatively fast growth of the phase variable $\varphi$. Conversely, the slowest expansion of the region $\Omega^{\rm T}(t)$ is observed in the controlled scenario where the local cell phenotype distribution remains concentrated around the sub-fittest phenotype $y=1.7$ (cf. brown curves, corresponding to the initial condition IC1 in Figure~\ref{fig::2}), as this corresponds to a slower growth of $\varphi$. Instead, an intermediate rate of expansion is observed in the other case (cf. orange curves, corresponding to the initial condition IC0 in Figure~\ref{fig::2}), since in this case the local cell phenotype distribution is initially centred at $\bar{y}_0=1.75$ but it then evolves into a Gaussian-like distribution centred at $y=1$ for all $x \in \Omega^{\rm T}(t)$ at later times, as shown by the results in the top panels of Figure~\ref{fig::1} (cyan curves) and in Figure~\ref{fig::3}. 

\begin{figure}
    \centering
    \includegraphics[width=0.32\textwidth, trim=0 5 45 15,clip]{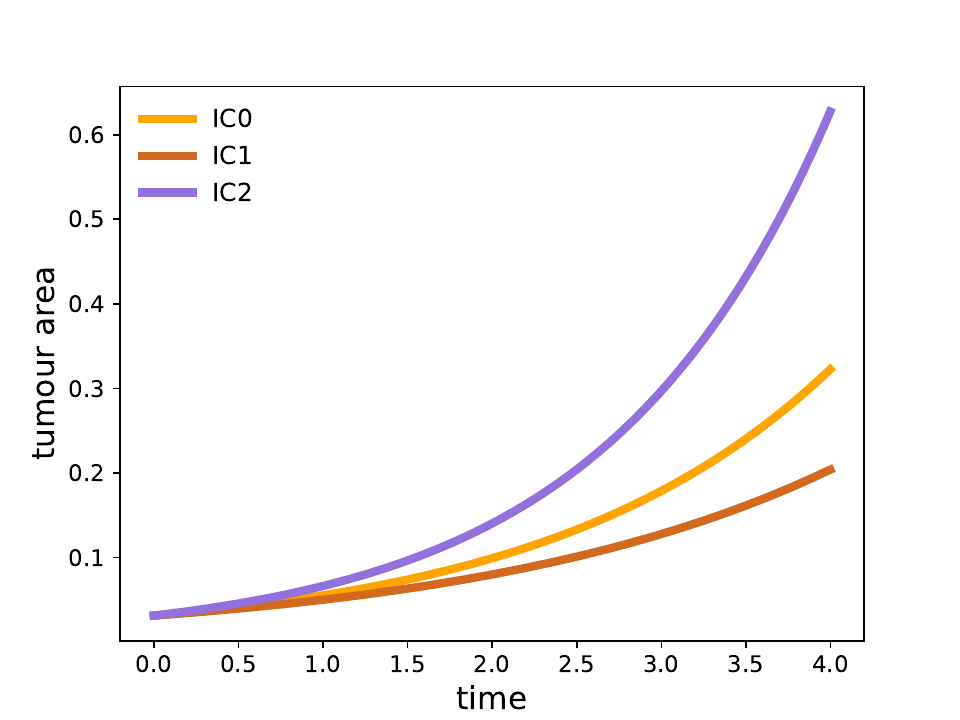}
	\includegraphics[width=0.32\textwidth, trim=0 5 45 15,clip]{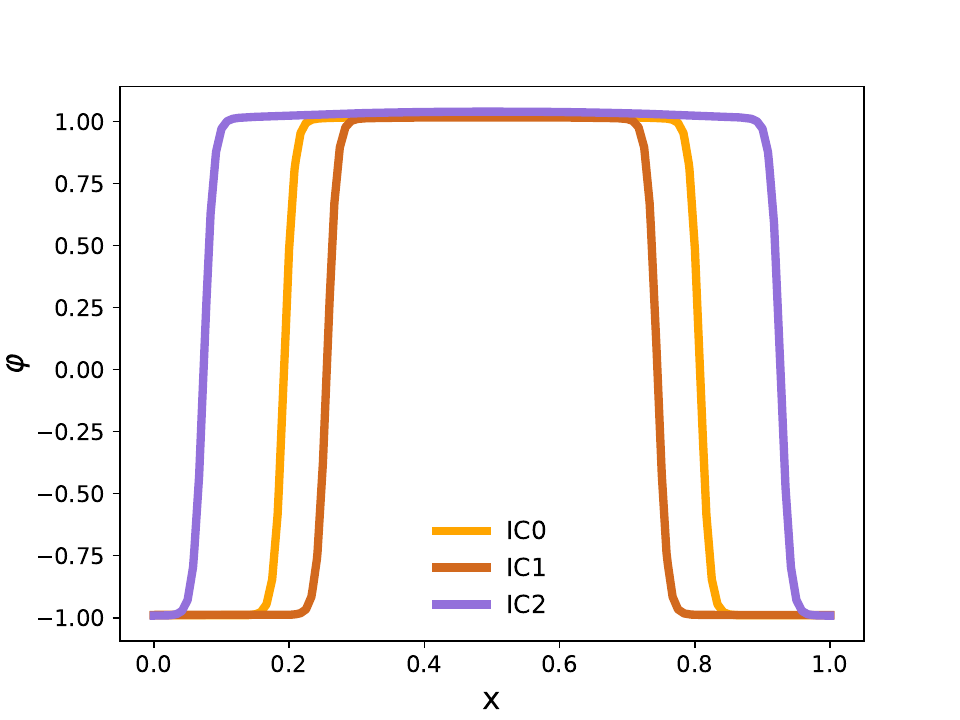}\\
	\hfil\includegraphics[width=0.6\textwidth]{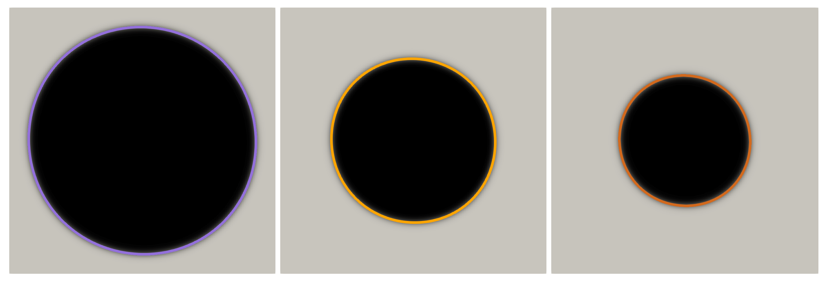}
	\includegraphics[width=0.1\textwidth]{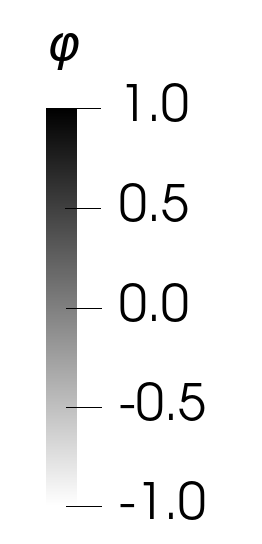}
    \caption{\small {\bf Top panels.} Plots of the measure of the set $\Omega^{\rm T}(t)$ defined via~\eqref{def:OmegaT} (left) and a cross section of the phase variable $\varphi$ at $t=4.0$ (right) when: $\alpha=0$ and $f_0$ is defined via~\eqref{def:f0} with $a=2.5$ and $\bar{y}_0=1$ (i.e. purple curves, for the initial condition IC2) or $\bar{y}_0=1.7$ (i.e. brown curves, for the initial condition IC1); $\alpha > 0$ (i.e. $\alpha= 5 \times 10^{2}$), for $\theta=0.5$ and $f_0$ defined via~\eqref{def:f0} with $a=2.5$ and $\bar{y}_0=1.75$ (i.e. orange curves, for the initial condition IC0). {\bf Bottom panels.}
	Plots of the phase variable $\varphi$ at $t=4.0$ for the initial condition IC2 (left), the initial condition IC0 (centre), and the initial condition IC1 (right).
	The superimposed coloured curves highlight the zero-isolines of $\varphi$.}
    \label{fig::2}
\end{figure}

Note that, in the controlled scenarios, $f(x,t,y) = f_0(y)$ for all $x \in \Omega^{\rm T}(t)$ and $t>0$, since $f_0$ is defined via~\eqref{def:f0} and $\alpha=0$, and thus $\displaystyle{\int_\YY \pp f(x,t,y) \, {\rm d}y =  \int_\YY \pp f_0(y) \, {\rm d}y}$, $\displaystyle{\int_\YY \qq f(x,t,y) \, {\rm d}y = \int_\YY \qq f_0(y) \, {\rm d}y}$, and $\displaystyle{\int_\YY \kk f(x,t,y) \, {\rm d}y = \int_\YY \kk f_0(y) \, {\rm d}y}$ for all $x \in \Omega^{\rm T}(t)$ and $t>0$. Hence, in these scenarios the phenotype-structured phase-field model~\eqref{eq:1}-\eqref{eq:4} can be regarded as reducing to a phase-field model of the form~\eqref{eq:01}-\eqref{eq:04} with $\displaystyle{p:=  \int_\YY \pp f_0(y) \, {\rm d}y}$, $\displaystyle{q := \int_\YY \qq f_0(y) \, {\rm d}y}$, and $\displaystyle{k:=\int_\YY \kk f_0(y) \, {\rm d}y}$. In contrast, in the case where $\alpha>0$, $f(x,t,y)$ evolves and, as a result of this, the values of $\displaystyle{\int_\YY \pp f(x,t,y) \, {\rm d}y}$, $\displaystyle{\int_\YY \qq f(x,t,y) \, {\rm d}y}$, and $\displaystyle{\int_\YY \kk f(x,t,y) \, {\rm d}y}$ evolve as well. This highlights a key difference of the the phenotype-structured phase-field model~\eqref{eq:1}-\eqref{eq:4} compared to the phase-field model~\eqref{eq:01}-\eqref{eq:04}, that is, its ability to capture dynamical changes in the model terms related to cell proliferation, death, and oxygen consumption due to the evolution of the phenotypic composition of the tumour.

\section{Conclusions and research perspectives}
\label{SEC:CONC}
In this paper, we have complemented the existing literature on phase-field models of tumour growth by incorporating inter-cellular phenotypic heterogeneity and the evolution of cell phenotypes into the phase-field modelling framework. This has been achieved by formulating a phenotype-structured phase-field model of nutrient-limited tumour growth. 

We have established existence and uniqueness of weak solutions for this model, within a functional framework that can accommodate a variety of biological situations. Moreover, we have presented a sample of numerical solutions intended to be a proof of concept for the ideas underlying the modelling framework. 

Future investigations may benefit from focused biological applications, thereby refining the model assumptions and allowing certain key questions to be addressed. To maintain simplicity while retaining generality, here we have assumed the rate of phenotype changes, $\theta$, to be constant. However, if this parameter was to be replaced by a phenotype-dependent function, in order to take into account the fact that cells with different phenotypes may undergo phenotype changes at different rates, the functional framework and the core estimates for the proof of the well-posedness result would remain intact. Furthermore, a natural extension of the present work would be to complement the numerical results presented here with analytical results on the qualitative properties of the solutions to the model equations. In this respect, appropriate extensions of the methods employed  in~\cite{Barlesetal2009,Bonnefonetal,Desv,lorenzipouchol2020} for analysing phenotype-structured integro-differential equations modelling evolutionary dynamics may prove useful to explore analytically key features of the local cell phenotype distribution $f$. In addition, it would be interesting to carry out pattern formation analysis to identify possible scenarios, corresponding to different assumptions on the model functions, under which the interplay between spatiotemporal and evolutionary mechanisms could induce breaking of the radial symmetry of $\Omega^{\rm T}(t)$, resulting in the formation of finger-like patterns, of the type of those produced by other partial differential equation models of tumour invasion into surrounding healthy tissue~\cite{Carteretal2025,Ciarlettaetal2011,Lietal2022,Lorenzietal2017}. It would also be interesting to investigate the existence of phenotype-structured travelling wave solutions, wherein the phase variable $\varphi$ behaves like a travelling front connecting $1$ to $-1$, while the local cell phenotype distribution $f$ takes the form of a (possibly sharp) Gaussian-like distribution, with the local mean phenotype that varies from place to place across the front. For this, we expect formal asymptotic methods building on those employed in~\cite{Lorenzietal2025b,Lorenzietal2022,LorenziPainter2022} to be useful.

The modelling approach proposed here to take into account inter-cellular variability in proliferation and death rates and in the rate of consumption of nutrients is flexible enough that it can possibly be adapted to incorporate phenotypic diversity of cells also in phase-field models of nutrient-limited tumour growth that take into account additional biophysical and biochemical mechanisms implicated in cancer invasion. For instance, it would be relevant to incorporate inter-cellular phenotypic diversity in phase-field models including chemotaxis effects. These are typically captured through a cross-diffusion term, as done for instance in~\cite{AS, CSS, EG, GLSS,  GSigSpr, KS2,  RSchS,SS}, which would lead to phenotype-structured phase-field models like the one defined by the following system:
\begin{alignat*}{2}
	& \dt \ph - \div \left(\mobm(\ph) \nabla \mu\right)
	=
	\hh(\ph) \left( \sigma \int_\YY \pp f \,{{\rm d}y} 
	-\int_\YY \qq f \,{{\rm d}y} \right)
	\qquad && \text{in $Q$,} \phantom{\int_\YY}
	\\
	& \mu = -\eps \Delta \ph +\frac1{\eps} F' (\ph) - \chi \sigma
	\qquad && \text{in $Q$,} \phantom{\int_\YY}
	\\
	& \dt f = \alpha\hh (\ph) \, \left[\theta \, \left(\int_{\cal Y} {\MM} f \,{{\rm d}y} - f\right) + {\left(\Rcel - \int_\YY \Rcel f \,{\rm d}y\right) f}\right]
	\qquad && \text{in $Q \times \YY$,} \phantom{\int_\YY}
	\\
	& \dt \sigma - D_\sigma\Delta \sigma
	+ \chi \Delta \ph
	+ \hh(\ph) \, \sigma \, \left(\int_\YY \kk f \,{{\rm d}y}\right)
	=
	b\left(\sigma_B - \sigma\right)
	\qquad && \text{in $Q$.}
\end{alignat*}
Boundary and initial conditions would then have to be prescribed accordingly. The terms $-\chi \sigma$ and $\chi \Delta \varphi$ in the above equations capture the tendency of tumour cells to move towards regions with higher nutrient concentrations.

Another class of models in which it would be interesting to incorporate a phenotypic structure would be phase-field models of nutrient-limited tumour growth wherein cells are transported by a velocity $\bv$ (i.e. the volume-averaged velocity of the mixture), which is prescribed by laws of fluid flow in porous media, such as Darcy's law~\cite{GLSS} or Brinkman's law~\cite{ CGSS1, EG,KS2}. This would lead to phenotype-structured phase-field models of the following type:
\begin{alignat*}{2}
	 & \div \bv = S_{\bv} \qquad && \text{in } Q,\\
	& \dt \ph + \div (\bv \ph) - \div \left(\mobm(\ph) \nabla \mu\right)
	=
	\hh(\ph) \left( \sigma \int_\YY \pp f \,{{\rm d}y} 
	-\int_\YY \qq f \,{{\rm d}y} \right)
	\,\, && \text{in $Q$,} \phantom{\int_\YY}
	\\
	& \mu = -\eps  \Delta \ph +\frac1 {\eps}  F' (\ph)-\chi \Delta\sigma
	\,\, && \text{in $Q$,} \phantom{\int_\YY}
	\\
	& \dt f = \alpha\hh (\ph) \, \left[\theta \, \left(\int_{\cal Y} {\MM} f \,{{\rm d}y} - f\right) + \left(\Rcel - \int_\YY \Rcel f \,{\rm d}y\right) f\right]
	\,\, && \text{in $Q \times \YY$,} \phantom{\int_\YY}
	\\
	& \dt \sigma + {\div (\bv \sigma)} - D_{\sigma} \, \Delta \sigma
	+ \chi \Delta \ph
	+ \hh(\ph) \sigma \left(\int_\YY \kk f \,{{\rm d}y}\right)
	=
	b\left(\sigma_B - \sigma\right)
	\,\, && \text{in $Q$.}	
\end{alignat*}
The velocity $\bv$ could then be given by the following relation, so as to place emphasis on factors like permeability, pressure drop, and flow rate,
\begin{alignat*}{2}
	& 
	\bv  = -  \nabla p  {-\ph\nabla (\mu + \chi \ph )}
	\qquad  && \text{in } Q,
\end{alignat*}
 or by Brinkman's law, that is,
\begin{alignat*}{2}
	& 
	- \div ( \eta (\ph) D\bv + \lambda(\ph)\div(\bv) \mathbf{I} ) + \nu(\ph) \bv + \nabla p = {-\ph\nabla (\mu + \chi \ph )}
	\qquad && \text{in } Q,
\end{alignat*}
which extends Darcy's law by introducing viscosity terms and is thus used in the presence of porous matrices and fluid-saturated region {(see, for instance, \cite{GLSS})}. Here, $S_{\bv}$ stands for a volumic source term arising from the mass balance of the system and it is related to the source terms on the \rhs\ of the equations for $\ph$ and $\sigma$. The term $D\bv = \frac 12 (\nabla \bv + \nabla \bv ^\top)$ is the symmetric gradient, $\mathbf{I}$ denotes the identity matrix, $p$ is the fluid pressure,  {and $-\ph\nabla (\mu + \chi \ph )$ captures nutrient contribution as well as the Korteweg force related to capillarity effects}. Moreover, $\eta(\ph)$ and $\lambda(\ph)$ are non-negative functions representing the shear viscosity and the bulk viscosity, respectively, while $\nu (\ph) $ stands for the permeability. This system would have again to be complemented with appropriate initial and boundary conditions.

Analysis and numerical simulation of phenotype-structured phase field models of tumour growth like~\eqref{eq:1}-\eqref{eq:4} and the above systems present significant challenges, which we put forward as potential avenues for future work promoting cross-fertilisation of research on Cahn--Hilliard-like equations and phenotype-structured equations, and acting as a catalyst of interdisciplinary research at the interface between biology and mathematics.

\section*{Acknowledgments}
TL gratefully acknowledges support from the Italian Ministry of University and Research (MUR) through the grant PRIN2022-PNRR project (No. P2022Z7ZAJ) ``A Unitary Mathematical Framework for Modelling Muscular Dystrophies'' (CUP: E53D23018070001) funded by the European Union - NextGenerationEU. GP and TL gratefully acknowledge support from the Istituto Nazionale di Alta Matematica (INdAM) and the Gruppo Nazionale per la Fisica Matematica (GNFM). GP also acknowledges support from the National Plan for Complementary Investments to the NRRP, project “D34H—Digital Driven Diagnostics, prognostics and therapeutics for sustainable Health care” (project code: PNC0000001), Spoke 4 funded by the Italian Ministry of University and Research. AS gratefully acknowledges partial support 
from the 
``MUR GRANT Dipartimento di Eccellenza'' 2023-2027,
from the GNAMPA (Gruppo Nazionale per l'Analisi Matematica, la Probabilit\`a e le loro Applicazioni) of INdAM (Isti\-tuto Nazionale di Alta Matematica)
project CUP E5324001950001, and from the Alexander von Humboldt Foundation.
 TL, GP, and AS gratefully acknowledge also the CNRS International Research Project `Mod\'elisation de la biom\'ecanique cellulaire et tissulaire' (MOCETIBI).


\vspace{3truemm}

\Begin{thebibliography}{10} \footnotesize

\bibitem{agosti2018personalized}
A. Agosti, C. Giverso, E. Faggiano, A. Stamm, P. Ciarletta,
A personalized mathematical tool for neuro-oncology: A clinical case study.
{\it Int. J. Non-Linear Mech.} {\bf 107} (2018), 170-181.

\bibitem{AS}
A. Agosti and A.~Signori,
Analysis of a multi-species Cahn–Hilliard–Keller–Segel tumor growth model with chemotaxis and angiogenesis.
{\it J. Differ. Equ.} {\bf 403} (2024), 308-367.

\bibitem{alnaes2015fenics}
M. Aln{\ae}s, J. Blechta, J. Hake, A. Johansson, B. Kehlet, A. Logg, C. Richardson, J. Ring, M.E. Rognes and G.N. Wells,
The FEniCS project version 1.5.
{\it Archive of numerical software} {\bf 3}(100) (2015)

\bibitem{Barlesetal2009}
G. Barles, S. Mirrahim and, B. Perthame,
Concentration in lotka-volterra parabolic or integral equations: a general convergence result. 
{\it Methods Appl. Anal.} {\bf 16} (2009) 321-340

\bibitem{Bonnefonetal}
O. Bonnefon, J. Coville and G. Legendre, 
Concentration phenomenon in some non-local equation. 
{\it Discrete Contin. Dyn. Syst. - B} {\bf 22} (2017), 763-781.

\bibitem{Burger2000}
R. B\"urger,
The mathematical theory of selection, recombination, and mutation. 
{\it  John Wiley \& Sons} (2000).

\bibitem{CH1}
J. W. Cahn,
On spinodal decomposition. 
{\it Acta Metall.} {\bf 9} (1961), 795-801.

\bibitem{CH2}
J. W. Cahn and J. E. Hilliard,
Free energy of a non-uniform system. 
{\it {I}. Interfacial free energy. J.
Chem. Phys.} {\bf 28} (1958), 258-267.

\bibitem{Carteretal2025}
P. Carter, A. Doelman, P. van Heijster, D. Levy, P. Maini, E. Okey and P. Yeung,
Deformations of acid-mediated invasive tumors in a model with Allee effect. 
{\it J. Math. Biol.} {\bf 90} (2025), 795-801.

\bibitem{chatelain2011morphological}
C. Chatelain, P. Ciarletta, and M.B. Amar,
Morphological changes in early melanoma development: influence of nutrients, growth inhibitors and cell-adhesion mechanisms.
{\it J. Theor. Biol.} {\bf 290} (2011), 46-59.

\bibitem{chen2014tumor}
Y. Chen and J.S. Lowengrub,
Tumor growth in complex, evolving microenvironmental geometries: a diffuse domain approach.
{\it J. Theor. Biol.} {\bf 361} (2014), 14-30.

\bibitem{Chisholmetal2016}
R.H. Chisholm, T. Lorenzi and J. Clairambault, 
Cell population heterogeneity and evolution towards drug resistance in cancer: biological and mathematical assessment, theoretical treatment optimisation. 
{\it Biochim. Biophys. Acta, Gen. Subj.} {\bf 1860} (2016), 2627-2645. 

\bibitem{Ciarlettaetal2011}
P. Ciarletta, L. Foret and M. Ben Amar, 
The radial growth phase of malignant melanoma: Muti-phase modelling, numerical simulation and linear stability. 
{\it J. R. Soc. Interface} {\bf 8} (2011), 345-368.

\bibitem{CGSS1}
P. Colli, G. Gilardi, A. Signori and J. Sprekels,
Cahn–Hilliard–Brinkman model for tumor growth with possibly singular potentials.
{\it Nonlinearity} {\bf 36} (2023), 4470-4500.

\bibitem{CSS}
P. Colli, A. Signori and J. Sprekels,
Optimal control of a phase field system modelling tumor growth with chemotaxis and singular potentials.
{\it Appl. Math. Optim.} {\bf 83} (2021), 2017-2049.

\bibitem{cristini2009nonlinear}
V. Cristini, X. Li, J.S. Lowengrub, S.M. Wise,
Nonlinear simulations of solid tumor growth using a mixture model: invasion and branching.
{\it J. Math. Biol.} {\bf 58} (2009), 723-763.

\bibitem{Desv}
L. Desvillettes, P.E. Jabin, S. Mischler and  G. Raoul,
On selection dynamics for continuous structured populations. 
{\it Commun. Math. Sci.} {\bf 6} (2008), 729-747.

\bibitem{EG}
M.~Ebenbeck and H.~Garcke,
Analysis of a Cahn--Hilliard--Brinkman model for tumour growth with chemotaxis.
{\it J. Differ. Equ.} {\bf 266} (2019),  5998-6036.

\bibitem{GARL_2}
H. Garcke and K.F. Lam,
Analysis of a Cahn--Hilliard system with non--zero Dirichlet 
conditions modeling tumor growth with chemotaxis.
{\it Discrete Contin. Dyn. Syst.} {\bf 37} (2017), 4277-4308.

\bibitem{GLS}
H.~Garcke, K.F.~Lam and A.~Signori,
On a phase field model of Cahn--Hilliard type for tumour growth with mechanical effects. 
{\it Nonlinear Anal. Real World Appl.} {\bf 57} (2021), 103192, https://doi.org/10.1016/j.nonrwa.2020.103192.

\bibitem{GLSS}
H.~Garcke, K.F.~Lam, E.~Sitka and V.~Styles,
A Cahn--Hilliard--Darcy model for tumour growth with chemotaxis and active transport.
{\it Math. Models Methods Appl. Sci.} {\bf 26} (2016),  1095-1148.

\bibitem{GSigSpr}
G. Gilardi, A. Signori and J. Sprekels,
Nutrient control for a viscous Cahn--Hilliard--Keller--Segel model with logistic source describing tumor growth.
{\it Discrete Contin. Dyn. Syst. Ser. S} {\bf 16} (2023), 3552-3572.

\bibitem{Hameletal2020}
F. Hamel, F. Lavigne, G. Martin and L. Roques,
Dynamics of adaptation in an anisotropic phenotype-fitness landscape.
{\it Nonlinear Anal. Real World Appl.} {\bf 54} (2020), 103107.

\bibitem{Huang2013}
S. Huang, 
Genetic and non-genetic instability in tumor progression: Link between the fitness landscape and the epigenetic landscape of cancer cells. 
{\it Cancer Metastasis Rev.} {\bf 32} (2013), 423-448.

\bibitem{KS2}
P. Knopf and A. Signori, 
Existence of weak solutions to multiphase Cahn--Hilliard--Darcy and Cahn--Hilliard--Brinkman models for stratified tumor growth with chemotaxis and general source terms. 
{\it Commun. Partial Differ. Equ.} {\bf 47} (2022), 233-278.

\bibitem{lady}
O.A. Lady\v zenskaja, V.A. Solonnikov and N.N. Uralceva, 
Linear and quasilinear equations of parabolic type.
{\it AMS Math. Monogr.} {\bf 23} (1968).

\bibitem{integrodiff}
V. Lakshmikantham,
Theory of integro-differential equations, {\it CRC press} {\bf 1} (1995).

\bibitem{Lietal2022}
H.-L. Li, S. Zhao and H.-W. Zuo,
Existence and nonlinear stability of steady-states to outflow problem for the full two-phase flow.
{\it J. Differ. Equ.} {\bf 309} (2022), 350-385.

\bibitem{logg2012dolfin}
A. Logg, G.N. Wells and J. Hake,
DOLFIN: A C++/Python finite element library.
{\it Automated Solution of Differential Equations by the Finite Element Method: The FEniCS Book}, (2012), 173-225.

\bibitem{Lorenzietal2025}
T. Lorenzi, K.J. Painter and C. Villa, 
Phenotype structuring in collective cell migration: a tutorial of mathematical models and methods. 
{\it J. Math. Biol.} {\bf 90} (2025).

\bibitem{Lorenzietal2025b}
 T. Lorenzi, F.R. Macfarlane and K.J. Painter,
 Derivation and travelling wave analysis of phenotype-structured haptotaxis models of cancer invasion.
{\it Eur. J. Appl. Math.} {\bf 36} (2025), 231-263.

\bibitem{Lorenzietal2022}
T. Lorenzi, B. Perthame and X. Ruan,
 Invasion fronts and adaptive dynamics in a model for the growth of cell populations with heterogeneous mobility.
{\it Eur. J. Appl. Math.} {\bf 33} (2022), 766-783.

\bibitem{LorenziPainter2022}
 T. Lorenzi and  K.J. Painter,
 Trade-offs between chemotaxis and proliferation shape the phenotypic structuring of invading waves.
{\it Int. J. Non Linear Mech.} {\bf 139}:103885 (2022).

\bibitem{Lorenzietal2017}
T. Lorenzi, A. Lorz and B. Perthame, 
On interfaces between cell populations with different mobilities. 
{\it Kinet. Relat. Models} {\bf 10} (2017), 299-311.

\bibitem{lorenzipouchol2020}
T. Lorenzi and C. Pouchol,
Asymptotic analysis of selection-mutation models in the presence of multiple fitness peaks.
{\it Nonlinearity} {\bf 33} (2020), 5791-5816.

\bibitem{Marusyketal2012}
A. Marusyk, V. Almendro and K. Polyak, 
Intra-tumour heterogeneity: A looking glass for cancer?
{\it Nature Rev. Cancer} {\bf 12} (2012), 323.

\bibitem{Michor2010}
F. Michor and K. Polyak, 
The origins and implications of intratumor heterogeneity. 
{\it Cancer Prevention Res.} {\bf 3} (2010), 1361-1364.

\bibitem{pozzi2022t}
G. Pozzi, B. Grammatica, L. Chaabane, M. Catucci, A. Mondino, P. Zunino and P. Ciarletta,
T cell therapy against cancer: A predictive diffuse-interface mathematical model informed by pre-clinical studies.
{\it J. Theor. Biol.} {\bf 547} (2022), 111172.

\bibitem{RSchS} 
E.~Rocca, G.~Schimperna and A.~Signori,
On a Cahn--Hilliard--Keller--Segel model with generalized logistic source describing tumor growth.
{\it J. Differ. Equ.} {\bf 343} (2023), 530-578.	

\bibitem{SS}
L.~Scarpa and A.~Signori,
On a class of non-local phase-field models for tumor growth with possibly singular potentials, chemotaxis, and active transport.
{\it Nonlinearity} {\bf 34} (2021), 3199-3250. 

\bibitem{S}
A. Signori,
Optimal distributed control of an extended model of tumor growth with logarithmic potential.
{\it Appl. Math. Optim.} {\bf 82} (2020), 517-549.

\bibitem{tierra2015numerical}
G. Tierra and F. Guill{\'e}n-Gonz{\'a}lez,
Numerical methods for solving the Cahn--Hilliard equation and its applicability to related energy-based models.
{\it Arch. Comput. Methods Eng.} {\bf 22} (2015), 269-289.

\bibitem{travasso2011phase}
R.D.M. Travasso, M. Castro and J. Oliveira,
The phase-field model in tumor growth.
{\it Philos. Mag.} {\bf 91} (2011), 183-206.

\bibitem{travasso2011phase}
R.D.M. Travasso, M. Castro and J. Oliveira,
The phase-field model in tumor growth.
{\it Philos. Mag.} {\bf 91} (2011), 183-206.

\bibitem{TsimringLevine1996}
L.S. Tsimring and H. Levine,
RNA virus evolution via a fitness-space model.
{\it Phys. Rev. Lett.} {\bf 76} (1996), 4440.

\bibitem{villaetal2021}
C. Villa, M.A.J. Chaplain and T. Lorenzi, 
Modelling the emergence of phenotypic heterogeneity in vascularised tumours. 
{\it SIAM J. Appl. Math.} {\bf 81} (2021), 434-453. 

\bibitem{villaetal2021b}
C. Villa, M.A.J. Chaplain and T. Lorenzi, 
Evolutionary dynamics in vascularised tumours under chemotherapy: Mathematical modelling, asymptotic analysis and numerical simulations. 
{\it Vietnam J. Math.} {\bf 49} (2021), 143-167. 

\bibitem{wise2008three}
S.M. Wise, J.S. Lowengrub, H.B. Frieboes and V. Cristini,
Three-dimensional multispecies nonlinear tumor growth -- I: model and numerical method.
{\it J. Theor. Biol.} {\bf 253} (2008), 524-543.

\End{thebibliography}

\End{document}
